\documentclass[11pt]{amsart}
\usepackage{amsmath,amssymb,amsthm}
\usepackage{mathtools}
\usepackage{graphicx}
\usepackage[margin=1in]{geometry}
\usepackage{hyperref}
\usepackage{enumerate}  % uncomment to use this package
\usepackage{upgreek,mathrsfs,comment}
\usepackage{color}

\newcommand{\ud}{\mathrm{d}}
\newcommand{\8}{\infty}

\newcommand{\bR}{\mathbb{R}}

\newcommand{\bE}{\mathbb{E}}
\newcommand{\bF}{\mathfrak{F}}
\newcommand{\bG}{\mathfrak{G}}

\newcommand{\sA}{\mathcal{A}}
\newcommand{\sB}{\mathcal{B}}
\newcommand{\sE}{\mathcal{E}}
\newcommand{\sF}{\mathcal{F}}
\newcommand{\sG}{\mathcal{G}}
\newcommand{\sH}{\mathcal{H}}
\newcommand{\sL}{\mathcal{L}}
\newcommand{\sM}{\mathcal{M}}
\newcommand{\sP}{\mathcal{P}}

\newcommand{\sT}{\mathcal{T}}

\newcommand{\blue}{\color{black}}
\newcommand{\black}{\color{black}}

\newcommand{\MP}{\textsf{MP}}
\DeclareMathOperator\ent{\upphi}
\DeclareMathOperator\Ent{\mathsf{Ent}}

\newcommand{\sfE}{\mathsf{E}}

\newcommand{\sfP}{\mathsf{P}}
\newcommand{\sfQ}{\mathsf{Q}}
\newcommand{\sfR}{\mathsf{R}}
\newcommand{\ddiv}{\overline{\rm div}\,}
\newcommand{\dnabla}{\overline{\nabla}}

\title[Optimal control for Markov Jump Processes]{Optimal control formulation of transition path problems for Markov Jump Processes}

\author{Yuan Gao}
\address{Department of Mathematics, Purdue University, West Lafayette, IN 47906}
\email{gao662@purdue.edu}

\author{Jian-Guo Liu}
\address{Department of Mathematics and Department of Physics, Duke University, Durham, NC 27708}
\email{jliu@math.duke.edu}

\author{Oliver Tse}
\address{Department of Mathematics and Computer Science, Eindhoven University of Technology, 5600~MB Eindhoven, The Netherlands}
\email{o.t.c.tse@tue.nl}



\date{}
\newtheorem{theorem}{Theorem}[section]
\newtheorem{proposition}[theorem]{Proposition}
\newtheorem{lemma}[theorem]{Lemma}

\theoremstyle{definition}
\newtheorem{definition}[theorem]{Definition}

\theoremstyle{remark}
\newtheorem{remark}[theorem]{Remark}
 \numberwithin{equation}{section}

\begin{document}

\begin{abstract}
   We formulate the transition path problem for Markov jump processes as a stochastic optimal control problem on path space. Transitions between metastable sets are induced by an unbounded terminal cost at a stopping time together with controlled modification of jump rates. The running cost takes an entropic form arising naturally from the Girsanov transform for jump processes, allowing both finite- and infinite-horizon formulations to be expressed as optimal changes of measure relative to a reference process.

We prove that the optimal path measure admits an explicit representation in terms of the committor function, which solves an elliptic boundary value problem. The singular control induced by the terminal cost is obtained via $\Gamma$-convergence, yielding a limiting controlled process whose transition rates correspond to a Doob-$h$ transform. The optimally controlled process generates transition paths almost surely while preserving the bridges of the reference process. 
\end{abstract}

\maketitle

\tableofcontents

\section{Introduction}

Rare events such as biochemical reactions, protein folding, genetic evolution, and quantum tunneling, despite their importance, happen with low probability. {\blue For example, nucleation is a rare event where Brownian molecules randomly form clusters, and once a critical size is reached, further binding triggers a phase transition \cite{zhang2016recent}. Other examples include transition paths in protein folding \cite{held2011mechanisms} where trajectories reach one region before another, and exit problems in chemical reactions where trajectories end in a region.  The low probability of such events makes both theoretical and computational studies on how quickly and through what mechanisms these rare events occur quite challenging. Specifically, calculating {\em transition paths}---possible paths taken during a transition---between two metastable sets, or simply sets of interest, has broad applications in physics, chemistry, and social science, and has thus become a key issue in applied mathematics \cite{weinan2006}.}

Many theoretical studies and computational techniques have been developed to tackle the transition path problem. Theoretically, this includes pathwise approaches based on large deviation theory \cite{olivieri-vares2005,freidlin-wentzell2012}, spectral theory \cite{davies-1-1982,davies-2-1982}, and potential-theoretic methods using capacities \cite{bovier-etal-2001,bovier-denHollander-2015}. On the computational side, approaches like the minimum action method \cite{MAM}, the (finite-temperature) string method \cite{string, FTSM}, and transition path theory (TPT) based on the \emph{committor function} \cite{weinan2006, MMS2009, Lu_Nolen_2015} are notable examples, each with its limitations.

A particularly successful approach is the committor function method within TPT, often used to determine the transition rate of paths. Using the committor function $h$, which solves a boundary value problem (see Section~\ref{sec:committor}), the Doob $h$-transform of an underlying stochastic process enables the study of a conditioned process used to predict transition paths connecting two metastable states. Although the committor function is widely used in computational transition path methods, the reasons and the ways in which it offers an efficient computational approach are not well-explained in the literature. Especially for rare events modeled by Markov jump processes, manipulations like adding an additional \emph{control field} or implementing selection procedures are necessary to facilitate (1) the efficient computation of transition paths, and (2) the generation of controlled paths that represent the most probable transition path for the original stochastic model, in line with large deviation theory similar to Freidlin-Wentzell or Schilder theory.

From a stochastic control perspective, the added control field should, in principle, optimize an appropriate running cost corresponding to the large-deviation rate function, along with a terminal cost, to induce the transition. For instance, \cite{GLLL23} formulated the transition path problem for an It\^o diffusion as a stochastic optimal control problem, with a closed-form solution where the optimal control field is explicitly expressed in terms of $\nabla \log h$. At a fixed noise intensity, the committor function provides an optimal drift that almost surely realizes the rare event. In this paper, we will make an analogous statement for Markov jump processes and examine the role of the committor function in various optimal control formulations.

\blue  
\subsubsection*{Journey from a proposed formulation to the main result} 
%This manuscript is written in a way that respects the chronological order in which our ideas developed. We made this choice to keep the other optimal control formulations---those not necessarily appropriate for the transition path problem---as they may be relevant to other problems. Moreover, it brings to light the thought process and motivation in the derivation of the final optimal control formulation for transition path problems.

To set the stage, we consider a Markov jump process $(X_t)_{t\ge0}$ on a Polish space $\varGamma$ with a time-homogeneous transition rate $L$, whose law is a probability measure $\sfR\in\sP(\Omega)$ on the Skorokhod space $\Omega\coloneqq D(\bR^+;\varGamma)$, $\bR^+\coloneqq[0,+\8)$ of $\varGamma$-valued left-limited and right-continuous (c\'adl\'ag) paths (cf.\ Section~\ref{sec:jump-process}). For any two nonempty, disjoint, measurable sets $A,B\subset\varGamma$, let $\tau_{AB}$ denote the first hitting time of the set $A\cup B$, i.e., $\tau_{AB} \coloneqq\inf\bigl\{t\geq 0 : \, X_t \in  A\cup B\bigr\}$. Further, set
	\begin{equation*}
    	f_{AB}:\varGamma\to[0,+\infty], \qquad f_{AB}(x) \coloneqq \begin{cases}
         +\infty & \text{for $x\in A$},  \\
         0 & \text{for $x\in B$}, \\
         \text{arbitrary} & \text{for $x\in (A\cup B)^c$}.
    \end{cases}
\end{equation*}
Since $\tau_{AB}$ forces $X_{\tau_{AB}}\in A\cup B$, only the values of $f_{AB}$ on $A\cup B$ affect the cost; the values on $(A\cup B)^c$ are irrelevant.

%{\red Recall diffusion processes and then lead to difficulties (including MC difficulties).}

As mentioned above, the transition path problem for diffusion processes has been extensively studied from both theoretical and computational viewpoints, primarily relying on the Freidlin-Wentzell large deviation principle. We will provide comparisons in later sections. In formulating an appropriate optimal control framework for the transition path problem in a Markov jump process, we encountered several additional challenges, which, though not unexpected, present significant difficulties:
\begin{enumerate}[(i)]
	\item Since drift and noise separation do not exist for Markov jump processes, unlike diffusion processes, we need to consider how to introduce a control variable for a given continuous-time Markov process and find a suitable cost function. Instead of the typical quadratic running cost, i.e., the kinetic energy of the control velocity field, we require an alternative running cost that measures the cost of changing the transition rate.
	\item A general Markov jump process that does not satisfy the detailed balance condition lacks a clear `energy landscape'. In contrast, for a reversible diffusion process, having an explicit global energy landscape enables interpreting the transition path as a geometric least action problem, where the action cost functional is based on a quadratic form inspired by the large deviation rate function for first exit time problems. The explicit energy landscape has also been used to derive a closed-form formula for an effective energy landscape that realizes transitions almost surely (cf.\ \cite{GLLL23}). Therefore, the second task is to propose a solvable stochastic optimal control problem for a general Markov jump process without detailed balance, where, in this context, solvability means that finding an optimal control can be reduced to solving a linear equation for the committor function. Moreover, the desired optimal control should not change the original transition paths (i.e., bridges) but should only stimulate their occurrence.
	\item Since the transition time of each trajectory is random, an appropriate optimal control problem should be considered over an infinite time horizon with a stopping time $\tau_{AB}$ as described above. This prevents formulating the control problem as a PDE-constrained deterministic optimal control problem for the time-marginal distribution of the process.
	\item The singular terminal cost due to the nature of the transition path problem introduces some technical difficulties. The sets $A, B$ are often chosen as metastable sets, where the reference process $X_t$ tends to visit with high probability. Therefore, to observe a transition from $A$ to $B$, the infinite penalty function $f_{AB}$ is added to prevent the event $X_{\tau_{AB}}\in A$. This infinite terminal cost $f_{AB}$ results in a singular control velocity whenever $X_t$ is near $A$, making it difficult to justify the existence of an optimally controlled process. A weak formulation (or martingale problem) for the controlled generator will be needed to establish the existence of an optimal controlled process.
\end{enumerate}

To resolve these issues, we begin our journey with a motivation for the quadratic cost function associated with It\^o diffusions via Girsanov's theory. Morally, if $(X_t)_{t\ge 0}$ is an It\^o diffusion with generator $\sL = \Delta/2 + b\cdot\nabla$ and law $\sfR$ and $(X_t^v)_{t\ge 0}$ is an It\^o diffusion with generator $\sL^v = \Delta/2 + (b+ v_t)\cdot\nabla$ and law $\sfP^v$ for some regular velocity field $v$, then Girsanov's theorem states that $\sfP^v$ has an $\sfR$-density such that (cf.~Section~\ref{sec:girsanov} or \cite{leonard2012girsanov})
\[
	\Ent(\sfP^v\,|\,\sfR) = \Ent(\sfP_0\,|\,\sfR_0) + \bE\biggl[\frac{1}{2}\int_0^\infty |v_t(X_t^v)|^2 \,\ud t \biggr].
\]
In particular, optimizing over admissible velocity fields $v$ w.r.t.\ to a quadratic cost function is equivalent to optimizing the entropy over admissible path measures $\sfP^v$ relative to $\sfR$. Hence, by considering the control problem on the level of path measures, we easily resolve difficulty (I) above, which, as a consequence, provides not only a candidate cost function but also the control mechanism (cf.\ Proposition~\ref{prop:girsanov}). 

Based on this discussion, and in view of the transition path  problem, we then propose the general stochastic optimal control problem
\begin{equation}\label{eq:soc-intro}
	\gamma_{\rm sto}(\mu)\coloneqq \inf \left\{\sfE_{\sfQ} \biggl[ f(X_\tau) + \log\frac{\ud \sfQ}{\ud \sfR} \biggr] \; :\; \sfQ\in \sP(\Omega),\;\sfQ_0=\mu \right\},\tag{{\sf s-OC$_\tau$}}
\end{equation}
which applies to both It\^o diffusions and jump processes for general terminal costs $f$ and stopping times $\tau$. An optimal path measure $\sfP\in\sP(\Omega)$ of \eqref{eq:soc-intro} serves as a solution to a martingale problem (cf.\ Definition~\ref{def:martingale}) with controlled transition rate $L_\sfP$. Therefore, in this setup, the search for a controlled process is equivalent to finding an optimal change of variables for the reference path measure $\sfR$.

Under relatively mild assumptions, we obtain our first main result (cf.\ Theorem~\ref{thm:soc}). Here, a function $f:\varGamma\to[0,+\infty]$ is called \emph{proper} if it is not identically $+\infty$, and its \emph{effective domain} is $\text{dom}\,f\coloneqq\{x\in\varGamma: f(x)<+\infty\}$.

\begin{theorem}\label{mainthm:soc}
	Let $(\Omega,\sF,\bF,\sfR)$ be a complete right-continuous filtered space and $\tau$ be an $\bF$-stopping time. Let $f:\varGamma\to[0,+\infty]$ be a proper $\sigma(X_\tau)$-measurable terminal cost satisfying additionally:
	\begin{align}
		\text{$\exists$ a $\sigma(X_\tau)$-measurable set $D\subset \text{dom}\,f$:}\quad \sfR_\tau(D)>0\quad\text{and}\quad \int_D f\,\ud\sfR_\tau <+\infty. \tag{$\mathsf{A}_f$}
	\end{align}
	Then, \eqref{eq:soc-intro} admits a minimizer for any $\mu\in\sP(\varGamma)$ with $\Ent(\mu\,|\,\sfR_0)<+\infty$.
\end{theorem}

Since $f\ge 0$, the finiteness of the cost yields the finiteness of the relative entropy for admissible path measures. In particular, Girsanov's theory offers a way to parameterize all such path measures using controlled transition rates $L^v = vL$ for some control velocity $v:\bR^+{\times} E_\varGamma\to[0,+\infty)$, $E_\varGamma\coloneqq \varGamma\times\varGamma$ (cf.\ \cite{leonard2012girsanov}). When $v$ is bounded and $(X_t^v)_{t\ge 0}$ is the associated controlled process with law $\sfP^v\in \sP(\Omega)$, then one can express the relative entropy as (cf.\ Proposition~\ref{prop:girsanov})
 \begin{align*}
  \Ent(\sfP^v\,|\,\sfR) = \Ent(\mu\,|\,\sfR_0) + \sfE_{\sfP^v}\biggl[\int_0^\tau\!\!\int_\varGamma\ent(v_t(X_{t^-}^v,y))\, L(X_{t^-}^v,\ud y)\,\ud t\,\biggr],
\end{align*}
where $\upphi(s)=s\log s-s+1$, which provides an explicit form for the running cost in terms of the velocity field $v$ and motivates the stochastic optimal control problem
\begin{equation}
	\begin{aligned}
		\inf_v\,\, &\sfE_{\sfP^v} \biggl[ f(X_{\tau}^v)+ \int_0^\tau\!\!\int_\varGamma\ent(v_t(X_{t^-}^v,y))\, L(X_{t^-}^v,\ud y)\,\ud t\biggr],\\
		&\text{s.t.} \quad \text{$X_t^v$\; has transition rate\; $L^v\coloneqq vL$},\qquad X_0^v\sim \mu .
	\end{aligned}\tag{{\sf ps-OC$_\infty$}}
\end{equation}

In Section \ref{sec:sto-tpt}, we finally deal with the transition path problem. More precisely, we focus on finding a closed representation for the minimizer $\sfP\!_{AB}$ of the variational problem \eqref{eq:soc-intro} with $f=f_{AB}$ and $\tau=\tau_{AB}$, and to determine the corresponding optimally controlled jump process via its transition rate.

Due to the singular nature of the terminal cost $f_{AB}$, the Girsanov transform presented in Proposition~\ref{prop:girsanov} cannot be applied directly and requires an approximation approach. The strategy involves first finding the optimal path measure $P^n_{AB}$ for a regularized problem using a cut-off terminal cost $f_{AB}^n$ and the disintegration formula. In this context, the Girsanov theorem can be applied to characterize the controlled transition rate $L_{AB}^n$ and, consequently, the associated control velocity field $v_{AB}^n$. By employing Gamma-convergence and noting that the sequence $(Z^n_\tau)_{n\ge 0}$ of $\sfR$-densities of $\sfP_{AB}^n$ is uniformly integrable for $n\ge 1$, we can let $n\to+\infty$ to deduce two consequences. 

Firstly, the sequence $(Z^n_\tau)_{n\ge 1}$ converges weakly to $Z_\tau$ in $L^1(\Omega,\sfR)$, where $\sfP\!_{AB}=Z_\tau\sfR$ is a minimizer of \eqref{eq:soc-intro} having finite relative entropy w.r.t.\ $\sfR$. In particular, $Z_\tau$ is a change of measure formula that takes the explicit form 
\[
	Z_\tau = \frac{h_{AB}(X_{\tau_{AB}})}{h_{AB}(X_0)},
\]
where the function $h_{AB}$ solves the boundary value problem \eqref{qn} below.

Secondly, the limit path measure $\sfP\!_{AB}$ is shown to solve the martingale problem with the controlled transition rate (cf.\ Theorem~\ref{thm:sto-rep})
\begin{equation*} 
	L_{AB}(x,\ud y) = \frac{h_{AB}(y)}{h_{AB}(x)} L(x,\ud y),\qquad x\in \varGamma,
\end{equation*}
from which we obtain the limit control velocity field $v_{AB}(x,y) = h_{AB}(y)/h_{AB}(x)$ used to generate a controlled process $X_t^{AB}$. This result can also be viewed as a generalization of the basic Girsanov transform in Proposition~\ref{prop:girsanov} for a class of jump processes with singular transition rates (cf.\ Remark~\ref{rem_sig}).

Together, these consequences imply our second main result.

\begin{theorem}\label{mainthm:tpt}
	Let the assumptions of Theorem~\ref{mainthm:soc} be satisfied with $f=f_{AB}$, $\tau=\tau_{AB}$ and $\text{supp}\, \mu\subset A^c$. Further, let $h_{AB}$ be the solution of \eqref{qn} and suppose
	\[
		\int_{\Gamma} \sfE_{\sfR^x}[e^{\tau_{AB}}]   \mu(\ud x) <+\infty,\quad%\text{for $\mu$-almost every $x\in \varGamma$},
	\]
	where $\sfE_{\sfR^x}[\cdot]\coloneqq \sfE_{\sfR}[\,\cdot\mid X_0=x]$ denotes the conditional expectation under $\sfR$ given the starting point $x\in\varGamma$.
	Then, \eqref{eq:soc-intro} admits a minimizer $\sfP\!_{AB}$ that takes the explicit form
		\[
			\sfP\!_{AB} = \int_\varGamma\frac{h_{AB}(X_{\tau_{AB}})}{h_{AB}(x)}\sfR^x\mu(\ud x)\in \sP(\Omega),
		\]
		and solves the martingale problem $\MP(\overline{L}_{AB},\mu)$ with L\'evy kernel
\begin{equation*}
	\overline{L}_{AB}(\omega,\ud t\,\ud y) \coloneqq \mathbf{1}_{[0,\tau_{AB})}(t)\frac{h_{AB}(y)}{h_{AB}(X_{t^-}(\omega))} L(X_{t^-}(\omega),\ud y)\,\ud t,\qquad \omega\in \Omega.
\end{equation*}

Consequently, the value function may be expressed as
\begin{equation*}
	\gamma_{\rm sto}(\mu) = \Ent(\mu\,|\,\sfR_0) -\int_\varGamma \log h_{AB}\,\ud \mu.
\end{equation*}
\end{theorem}

We emphasize that the optimal path measure $\sfP\!_{AB}$ does not alter the original behavior of the transition paths. More precisely, the disintegration representation for the optimal path measure $P_{AB}$ indicates that the bridges of $P_{AB}$ coincide with the bridges of the reference measure $\sfR$. This result is parallel to transition paths for drift-diffusion processes found in \cite{Lu_Nolen_2015}, and it is similar to the Schr\"odinger bridge problem \cite{leonard-2012}.
Our result also provides a theoretical guarantee for using the committor function $h_{AB}$ to control jump processes in the computations of transition paths. We note that solving the committor equation for high-dimensional problems is computationally expensive. However, in the context of importance sampling, our variational perspective shows that it is not necessary to solve the committor equation exactly. Even an approximate solution can significantly improve the computation of transition paths. This can be seen as a form of preconditioning, where the importance sampling strategy offers effective guidance for exploring transition pathways \cite{BD_book}. \black

\subsubsection*{Comparison to other work}

We also discuss some other related methodologies in the study of rare events and stochastic optimal control problems. As a classical method, the exponential change of variables can be used to connect nonlinear Hamilton-Jacobi equation (HJE) and linear equations, known as the Cole-Hopf transformation in differential equations. The Hamilton-Jacobi method is also linked to the Lagrangian perspective in the form of an optimal control problem through variational formulas \cite{Sheu85, fleming06, GL22}, which is also called the prelimit variational representation \cite{BD_book}. \blue Particularly, the logarithmic transform was adopted in \cite{Sheu85} (see also \cite{fleming06}) to derive the HJE and study the asymptotic behaviors for exit probability of jump processes on $\bR^d$. The optimal control formula with a Doob-transformation type controlled generator $L^h f\coloneqq h^{-1}L(hf)$ was derived. In contrast, our work generalizes the exit problem to transition path problem (where trajectories reach one region before another) for jump processes on general Polish spaces $\varGamma$. We linked stopping-time problems to variational problems via Girsanov theory, meanwhile we justified  the associated singular optimal control due to the singular terminal cost at one region. \black  As noted in \cite{BD_book}, the relative entropy is widely used in selecting an appropriate cost structure (Lagrangian) in the variational representations. These representations, after taking a zero noise limit, are widely used in the study of large deviation principles \cite{varadhan1984large}, particularly for chemical reactions modeled by Markov jump processes \cite{Kurtz15, Dembo18, GL22t}. A recent related result in transition path theory for drift-diffusion processes is \cite{Lu_Nolen_2015}, which established the uniqueness and existence of a conditioned process (an optimally controlled process possibly starting from $\partial A \cup \partial B$) with a singular generator $L^h f\coloneqq h^{-1}L(hf)$ in the form of a strong solution to a Brownian motion-driven SDE up to the stopping time $\tau_{AB}$. Our Theorem \ref{thm:sto-rep} for pure jump processes can be viewed as a parallel result using a martingale problem in a weak formulation with a similar singular generator (L\'evy kernel)  $\overline{L}_{AB}$. From a numerical perspective, the controlled Markov chain via the Doob transformation has already been used in previous algorithms for transition path computations; see \cite{Hartmann2016, MMS2009, GLLL23} for computational methods concerning high-dimensional committor functions.

\subsubsection*{Outline} \blue
 The remaining contents are organized as follows. Section~\ref{sec2} provides an overview of the preliminaries related to Markov jump processes, the martingale problem, and the transition path problem, emphasizing the concept of the committor function.
 Section~\ref{sec:girsanov} introduces the Girsanov transformation for pure jump processes with bounded control velocity, motivated by the optimal change of measures for diffusion processes. This section establishes a connection between the running cost and a relative entropy for path measures.
 In Section~\ref{sec:sto-opt}, we present the stochastic optimal control formulation for the transition path over an infinite time horizon. The running cost is defined by the relative entropy for path measures, with a (possibly) unbounded terminal cost at a stopping time. We demonstrate the existence of an optimal path measure (Theorem \ref{thm:soc}).
 Section~\ref{sec:sto-tpt} returns to the specific transition path problem, offering a closed formula for the optimal path measure derived from the committor function. We also derive the associated martingale problem with a singular generator (see Theorem~\ref{thm:sto-rep}).
 Section~\ref{sec5} concludes the document with final remarks and observations.
 Lastly, Appendix~\ref{app:girsanov} contains a self-contained proof of the Girsanov transform as detailed in Proposition~\ref{prop:girsanov}, while Appendix~\ref{sec:deterministic} discusses the optimal control formulation in a finite time horizon, along with its closed-form solution.
 \black

\section{Preliminaries}\label{sec2}
This section describes transition path problems for continuous time Markov jump processes, the committor function, and its probabilistic representation.

\subsection{Markov jump processes, local attractors, and transition paths}\label{sec:jump-process}
This subsection reviews and outlines the basic setup for the transition path problem, focusing primarily on Markov jump processes with bounded transition rates.

\medskip

\paragraph{\bf \em Jump processes on $\varGamma$} 
	Let $\varGamma$ be a Polish space and consider the Skorokhod space $\Omega\coloneqq D(\bR^+;\varGamma)$ of all c\`adl\`ag paths from $\bR^+$ to $\varGamma$. The canonical (or coordinate) process $X=(X_t)_{t\geq 0}$ is defined by
\[
	 \bR^+ \times \Omega\ni (t,\omega)\mapsto X_t(\omega) = \omega_t \in\varGamma,
\]
i.e., $X$ is the identity map on $\Omega$ and $X_t:\Omega\to \varGamma$ is its $t$-marginal. The space $\Omega$ is equipped with the $\sigma$-algebra $\sF^X\coloneqq\sigma(X_t:\,t\ge 0)$ generated by the $t$-marginals $X_t$, $t\ge 0$, and the canonical filtration $\bF^X=\{\sF_t^X\coloneqq\sigma(X_s:s\le t)\}_{t\ge 0}$. In view of hitting times, we also consider the right-continuous extension $\bF_+^X$ of the canonical filtration $\bF^X$. For a given path measure $\sfP\in\sP(\Omega)$, we can complete the canonical $\sigma$-algebra $\sF^X$ and the right-continuous filtration $\bF^X_+$ to obtain a complete right-continuous filtered space $(\Omega,\sF,\bF,\sfP)$, where $\sF$ and $\bF$ are the $\sfP$-completion of $\sF^X$ and $\bF_+^X$ respectively. We further denote by $\sT\coloneqq\{t\ge 0: X_t\ne X_{t^-}\}$ the set of canonical jump times and $\sT_t\coloneqq \sT\cap [0,t]$ the jump times restricted to the interval $[0,t]$ for each $t\ge 0$, with the convention $\sT_0=\emptyset$. 

\medskip

An $\bF$-predictable \emph{Markovian L\'evy kernel} $\overline L$ is a random nonnegative measure on $\bR^+\times \varGamma$ taking the form
\[
	\overline{L}(\omega,\ud t\,\ud y)= \,L_t(X_{t^-}(\omega),\ud y)\,\ud t,\qquad \omega\in \Omega,
\]
such that the filtration $\bG=\{\sG_t\}_{t\ge 0}$ defined by
\[
	\sG_t \coloneqq \sigma\bigl(\,\overline L((0,s]\times A):s\le t,\; A\subset\varGamma\;\;\text{Borel}\bigr),
\]
satisfies $\bG\subset\bF$, and $\{L_t(x,\cdot)\}_{(t,x)\in \bR^+\times\varGamma}$ is a Borel family of Radon measures on $\varGamma$.

\medskip

Adapted from \cite{Jacod1975} (see also \cite{Kallenberg2021,leonard2012girsanov}), we make the following definition.

\begin{definition}[Martingale problem]\label{def:martingale}
	A path measure $\sfP\in\sP(\Omega)$ is said to solve the martingale problem $\MP(\overline{L},\mu)$ with Markovian L\'evy kernel $\overline L$ if there exists a right-continuous filtration $\bF$ such that $\overline L$ is $\bF$-predictable and the following properties hold:
	\begin{enumerate}
		\item $\sfP$ has initial law $\mu\in \sP(\varGamma)$, i.e., $\sfP(X_0\in \cdot)=\mu$, and
		\item for every bounded Borel function $\varphi\in B_b(\bR^+\times E_\varGamma)$, the $\bR$-valued process
	\[
		\sum_{s\in\sT_t} \varphi_s(X_{s^-},X_s) - \iint_{(0,t]\times \varGamma} \varphi_s(X_{s^-},y)\, \overline{L}(\ud s\,\ud y),
	\]
is a $\sfP$-martingale w.r.t.\ $\bF$. 
	\end{enumerate}
	In this case, we say that $\sfP$ admits the L\'evy kernel $\overline{L}$.
\end{definition}

\medskip

Let $L:\varGamma\times \sM_+(\varGamma)\to[0,+\infty)$ be a bounded \emph{transition rate} satisfying
\begin{equation}\label{eq:ref-levy1}
 c_L\coloneqq\sup_{x\in\varGamma} \int_{\varGamma} L(x,\ud y) <+\infty.  
\end{equation}	
Then for the Markovian L\'evy kernel 
\begin{align}\label{eq:ref-levy}
	\overline{L}_\sfR(\omega,\ud t\,\ud y) \coloneqq L(X_{t^-}(\omega),\ud y)\,\ud t,\qquad \omega\in \Omega,
\end{align}
there is a unique {\em reference measure} $\sfR\in\sP(\Omega)$  admitting the L\'evy kernel $\overline{L}_\sfR$ \cite[Theorem~3.6]{Jacod1975}. We further assume that the canonical process $X$ under $\sfR$ is {\em positive} (in the sense of \cite[Section~10]
{Tweedie2009}), which therefore implies the existence of a unique invariant measure $\pi\in\sP(\varGamma)$ for $\sfR$; see \cite[Section~10]{Tweedie2009}. When $\varGamma$ is a finite or countable space, then recurrence and {\em strong aperiodicity} suffice to guarantee the existence of a unique invariant measure for $\sfR$ \cite[Proposition 10.4.2]{Tweedie2009}. A set $M$ is called $\sfR$-{\em metastable} if there exists a probability measure $\pi_M$ satisfying
\[
	\lim_{t\to\infty} \sfR^x\bigl(X_t\in A\;|\;\sigma_M>t\,\bigr) = \pi_M(A),\qquad \forall\,x\in M,\;A\in \sB(\varGamma),
\]
where $\sigma_M$ is the first return time of the process $X$ to the set $M$ \cite{Tweedie2009}. Here, $\sfR^x$ denotes the conditional probability for paths starting at $x\in\varGamma$, i.e., $\sfR^x=\sfR(\cdot\,|\,X_0=x)$.

For the continuous case, the overdamped Langevin process is described by the stochastic differential equation  $\ud X_t = -\nabla U(X_t)\,\ud t + \ud B_t$, where the noise is given by a Brownian motion and the drift is in the gradient form $-\nabla U$. In this case, the wells around the local minima of the energy landscape $U$ define potential metastable sets. For both continuous and jump processes, these regions describe typical states that occur, e.g., in chemical reactions. The transition path problem aims to study the transitions between two metastable sets, which can be regarded as rare events. 

As mentioned in the introduction, for general Markov jump processes, \blue one does not have a geometric least action problem with a quadratic Lagrangian cost function. Unlike the diffusion case, where control appears as an additional velocity field whose kinetic energy needs to be minimized, the formulation of the running cost for jump processes is unclear. In particular, one needs to propose an appropriate running cost that does not alter the bridges when regulating the transition paths. Moreover, for diffusion processes, the corresponding optimal control problem with a quadratic Lagrangian is well-studied and admits explicit closed-form solutions. \black Therefore, it is both interesting and important to propose a solvable optimization problem that facilitates transitions between multiple metastable sets. The committor function plays a crucial role in this context, as it describes the probability of the Markov process hitting one metastable set before another. We will explore a stochastic optimal control formulation to realize the transition path with minimal running cost, where the cost is given by the relative entropy, allowing us to reduce the optimization problem to a linear one. The committor function will be reviewed in detail in the next section.

\subsection{The committor function}\label{sec:committor}
%In this subsection, we will precisely introduce and discuss the committor function.

Consider two nonempty, disjoint subsets $A,B\subset\varGamma$. In practice, these sets would be metastable sets for $\sfR$.
If $\varGamma$ is countable, we can assume without loss of generality that $A$ and $B$ are simply two distinct points in $\varGamma$. The aim is to compute the most probable transition path from $A$ to $B$, such that it reaches $B$ before returning to $A$. 

The {\em committor function} $h_{AB}$ (or \emph{equilibrium potential} in the potential theory community)  of the pair $(A,B)$ is the unique solution of the following boundary value problem
\begin{equation}\label{qn}
	\begin{cases}\displaystyle
		\;\; \sL h(x)\coloneqq\int_{\varGamma} \dnabla h(x,y)\,L(x,\ud y) = 0 & \quad x\in  (A\cup B)^c, \\
		\hspace*{7.4em} h(x)=\mathbf{1}_B(x) & \quad x\in A\cup B,\tag{{\sf BVP}}
	\end{cases}
\end{equation}
where $\dnabla h(x,y)=h(y)-h(x)$, $x,y\in\varGamma$ denotes the discrete gradient, \blue $\sL$ is the generator corresponding to the reference process $\sfR$. \black  The unique solvability of \eqref{qn} for general transition rates on general sets $A,B$ is a non-trivial problem (cf.\ \cite{bovier-denHollander-2015} for some results). However, since it is not the main focus of this work, we henceforth assume that a unique solution of \eqref{qn} exists for any measurable pair $(A,B)$.

Due to Dynkin's formula \cite[Theorem 9.21]{Fima_book}, the committor function can be given a probabilistic interpretation. Indeed, defining the stopping time
\[\tau_{AB} \coloneqq\inf\bigl\{t\geq 0 : \, X_t \in  A\cup B\bigr\},\]
and assuming that $\sfE_{\sfR^x}[\tau_{AB}]<+\infty$, Dynkin's formula then gives for any $x\in (A\cup B)^c$,
\begin{equation}\label{eq:dynkin}
    \sfE_{\sfR^x}\bigl[h_{AB}(X_{\tau_{AB}})\bigr]- h_{AB}(x) = \sfE_{\sfR^x} \biggl[\int_0^{\tau_{AB}} \!\!\!\int_{\varGamma} \,\dnabla h_{AB}(X_{t^-},y)\,\overline{L}_\sfR(\ud s\,\ud y)\biggr] =0.
\end{equation}
Thus, the committor function provides the probability of hitting $B$ before $A$, i.e.
\begin{equation}\label{h_prop}
    h_{AB}(x) = \sfE_{\sfR^x}\bigl[\mathbf{1}_B(X_{\tau_{AB}})\bigr]=\sfR^x(\tau_B<\tau_A),\qquad x\in  {(A\cup B)}^c,
\end{equation}
where $\tau_A$ and $\tau_B$ are the hitting times of the sets $A$ and $B$ respectively.

The reaction rate, the probability density, and the reactive current (probability flux) of transition paths can be computed based on the committor function. For diffusion processes, these formulas were first obtained in \cite{weinan2006} and then theoretically studied in \cite{Lu_Nolen_2015}. Particularly, the process conditioned on hitting $B$ before $A$ has a distribution $\sfP^x$ on $\Omega$ that solves a martingale problem starting from $x$ with the (infinitesimal) generator $\sL^hf \coloneqq h^{-1}L(hf)$ given via the $h$-transform \cite[Section~7.2]{Pinsky} (see also \cite[Section~4.1]{Pinsky} for an explicit form of the conditioned transition probabilities). Here $L(hf)(x)\coloneqq \int_\varGamma(hf(y)-hf(x))L(x,\ud y)$ uses the transition rate $L$, and $\sL^h$ is to be distinguished from $L$. For the transition path problem of diffusion processes, \cite[Theorem 1.1]{Lu_Nolen_2015} showed that the conditioned process (possibly starting from $\partial A \cup \partial B$) with singular generator $\sL^h$ still has a unique strong solution up to the stopping time $\tau_{AB}$. Therefore, with the knowledge of the committor function, one can construct many numerical stochastic algorithms for the transition path problem; see \cite{Hartmann2016,  GLLL23} for efficient Monte Carlo simulations using controlled processes and \cite{MMS2009} for the dominated transition path method by finding bottleneck of transition capacity recursively.

\section{Stochastic Optimal Control Formulation}\label{sec:sto-opt}

\blue In this section, we formulate the original transition path problem for jump processes as a stochastic optimal control problem based on path measures. First, we draw inspiration from the stochastic optimal control formulation for transition path problems driven by It\^o diffusions. Using the Girsanov transformation for It\^o diffusions, we connect the kinetic-energy-type running cost for the control velocity to the relative entropy of path measures in Section~\ref{sec:girsanov}. After introducing the Girsanov transformation for jump processes, we define the optimal control problems for jump processes using the relative entropy of path measures on $\Omega$, which introduces an entropy-based running cost for the control velocity associated with a {\em controlled} L\'evy kernel. Second, the change of measure problem can be further simplified to convex optimization problems over path measures with various constraints, such as those with a fixed initial condition or fixed terminal cost at a stopping time. After accounting for a possibly unbounded terminal cost, we present a general existence and uniqueness result in Theorem~\ref{thm:soc} for the optimization problem. \black

In what follows, the relative entropy will play a central role. For a measurable space $\Upsilon$ and two measures $\alpha,\beta\in \sM_+(\Upsilon)$, the relative entropy of $\alpha$ w.r.t.\ $\beta$ is defined as
\begin{align}\label{eq:def-entropy}
	\Ent(\alpha\,|\,\beta)\coloneqq \begin{cases}\displaystyle
		\int_{\Upsilon}\, \ent\left(\frac{\ud \alpha}{\ud \beta}\right) \ud \beta & \text{if $\alpha\ll \beta$}, \\
		+\infty &\text{otherwise},
	\end{cases}\qquad \ent(s) =s\log s -s +1,
\end{align}
where $\alpha\ll \beta$ means that $\alpha(A)=0$ whenever $\beta(A)=0$, for any measurable set $A\subset\Upsilon$.

Moreover, for any measure $\sfP\in\sP(\Omega)$ and measurable map $\phi:\Omega\to \Omega_\phi$, we set $\sfP_\phi\coloneqq\phi_\#\sfP\in\sP(\Omega_\phi)$ and define $\sfP^\phi$ to be the conditional measure obtained via the disintegration theorem, i.e.
\begin{align}\label{def:disintegration}
	\sfP(A) = \int_{\Omega_\phi} \sfP^{\phi=\eta}(A)\,\sfP_\phi(\ud \eta)\qquad\text{$\forall$ measurable sets $A\subset\Omega$}.
\end{align}
We note that if $\sfP,\sfR\in\sP(\Omega)$ with $\sfP\ll \sfR$, then the disintegration theorem also gives \cite{leonard2014path}
\[
	\frac{\ud \sfP}{\ud \sfR}(\omega) = \frac{\ud \sfP_\phi}{\ud \sfR_\phi}(\phi(\omega))\,\frac{\ud \sfP^{\phi}}{\ud \sfR^{\phi}}(\omega)\qquad\text{for $\sfP$-almost every $\omega\in\Omega$.}
\]
In the following, we use the shorthand notation $\sfP = \sfP_\phi\otimes \sfP^{\phi}$ to denote \eqref{def:disintegration}.

We further mention a fundamental additive property of the relative entropy under disintegration:
\begin{align}\label{eq:entropy-disintegrate}
	\Ent(\sfP\,|\,\sfR) = \Ent(\sfP_\phi\,|\,\sfR_\phi) + \int_{\Omega_\phi} \Ent(\sfP^{\phi=\eta}\,|\,\sfR^{\phi=\eta})\,\sfP_\phi(\ud \eta).
\end{align}

\subsection{Girsanov's theorem for Markov jump processes}\label{sec:girsanov}

\blue 
We first motivate the design of the running cost, as derived from Girsanov's theorem, for diffusion processes, cf. \cite{leonard2012girsanov}.  
\begin{proposition}[Girsanov's theorem for diffusion processes with bounded velocities]
 Let $(\Omega,\sF,\bF,\sfR)$ be a complete  continuous filtered space and $v:\bR^{+}\times  \bR^d\to \bR^d$ be a velocity field satisfying
\[
    |v_t(x)|\le \overline{v}<+\infty\qquad\forall\,(t,x)\in \bR^+\times  \bR^d.
\]
Suppose $\sfR$ solves the martingale problem with generator $\sL= \Delta/2 + b\cdot \nabla$, where $b:\bR^d\to \bR^d$ is a regular drift field.
Then, the $\bR$-valued process $Z^v=(Z_t^v)_{t\ge 0}$ defined by
\[
	Z_t^v \coloneqq \exp\Biggl( \int_0^t v_s(X_s) \cdot \ud B_s -   \frac{1}{2}\int_0^t |v_s(X_s)|^2\, \ud s\Biggr)  ,
\]
is a non-negative $\sfR$-martingale w.r.t.\ $\bF$. 

Moreover, the path measure $\sfP^v\in\sP(\Omega)$ with initial law $\mu\in L^1(\varGamma,\sfR_0)$ defined by
\begin{align}\label{eq:P-girsanov-formula-diff}
	\sfP^v \coloneqq \frac{\ud\mu}{\ud \sfR_0}(X_0)\,Z_t^v\,\sfR\quad\text{on $\sF_t$,\;\;$t\ge 0$},
\end{align}
solves the martingale problem $\MP(\overline{L}^v,\mu)$ with generator
$\sL^v = \Delta/2 + b \cdot \nabla + v_t \cdot \nabla.$

We also know $Z_\infty^v\coloneqq \lim_{t\to+\infty} Z_t^v$ exists in $L^1(\Omega,\sfR)$,
\begin{align}\label{eq:entropy-formula-diff}
	\Ent(\sfP^v\,|\,\sfR) = \Ent(\mu\,|\,\sfR_0) + \sfE_{\sfP^v}\biggl[\frac{1}{2}\int_0^\infty |v_s(X_s)|^2\, \ud s\,\biggr],
\end{align}
and $\sfE_\sfR[Z_\infty^v\,|\,\sF_\tau]=Z_\tau^v$ for any $\bF$-stopping time $\tau$.
\end{proposition}

From this Girsanov theorem, it is natural to design the running cost as a quadratic form (kinetic energy) of the additional drift velocity $v$, which controls the original diffusion process. It also indicates that the commonly used optimal control formulation  of the transition path problem for diffusion processes   (cf. \cite{Lu_Nolen_2015, Hartmann2016,    GLLL23})
\begin{align*}
    \inf_{v}\,\, &\bE \biggl[f(X_\tau^v) +  \frac{1}{2}\int_0^\tau |v_s(X_s^v)|^2\, \ud s\,\biggr],\\[0.5em]
    &\text{s.t.} \quad \ud X_t^v = v_t(X_t^v)\,\ud t + b(X_t^v)\,\ud t + \ud B_t, \qquad X_0^v\sim \mu,
\end{align*}
is equivalent to an optimal change of path measure formulation, where the relative entropy between path measures is minimized.

Motivated by this, we revisit the Girsanov theorem for jump processes with a bounded control velocity, and then introduce the stochastic optimal control formulation for jump processes in terms of both path measures and the controlled L\'evy kernel. 
\black

Our first result is the following statement, which is known in various forms in the literature. We include its proof in Appendix~\ref{app:girsanov} for completeness.

\begin{proposition}[Girsanov's theorem for bounded velocities]\label{prop:girsanov} Let $(\Omega,\sF,\bF,\sfR)$ be a complete right-continuous filtered space and $v:\bR^+\times E_\varGamma\to[0,+\infty)$ be a Borel-measurable velocity field satisfying
\[
	0<\underline{v}\le v_t(e)\le \overline{v}<+\infty\qquad\forall\,(t,e)\in \bR^+\times E_\varGamma.
\]
Then, the $\bR$-valued process $Z^v=(Z_t^v)_{t\ge 0}$ defined by
\[
	Z_t^v \coloneqq \exp\left(- \iint_{(0,t]\times \varGamma} \bigl(v_s(X_{s^-},y)-1\bigr)\,\overline L_\sfR(\ud s\,\ud y)\right)\prod_{s\in \sT_{t}} v_s(X_{s^-},X_s) ,
\]
is a non-negative $\sfR$-martingale w.r.t.\ $\bF$.

Moreover, the path measure $\sfP^v\in\sP(\Omega)$ with initial law $\mu\in L^1(\varGamma,\sfR_0)$ defined by
\begin{align}\label{eq:P-girsanov-formula}
	\sfP^v \coloneqq \frac{\ud\mu}{\ud \sfR_0}(X_0)\,Z_t^v\,\sfR\quad\text{on $\sF_t$,\;\;$t\ge 0$},
\end{align}
solves the martingale problem $\MP(\overline{L}^v,\mu)$ with the L\'evy kernel 
\[
	\overline{L}^v(\omega,\ud t\,\ud y) \coloneqq L_t^v(X_{t^-}(\omega),\ud y)\,\ud t,\quad \omega\in \Omega,
\]
where $L_t^v(x,\ud y)\coloneqq v_t(x,y)L(x,\ud y)$.

 If in addition $\Ent(\sfP^v\,|\,\sfR)<+\infty$, then $Z_\infty^v\coloneqq \lim_{t\to+\infty} Z_t^v$ exists in $L^1(\Omega,\sfR)$,
\begin{align}\label{eq:entropy-formula}
	\Ent(\sfP^v\,|\,\sfR) &= \Ent(\mu\,|\,\sfR_0) + \sfE_{\sfP^v}\biggl[\iint_{\bR^+\times \varGamma}\ent\biggl(\frac{\ud \overline{L}^v}{\ud \overline{L}_\sfR}\biggr)\, \ud \overline{L}_\sfR\,\biggr]\\
    &= \Ent(\mu\,|\,\sfR_0) +\sfE_{\sfP^v}\biggl[ \iint_{\bR^+\times \varGamma} \ent(v_t(X_{t^-},y)) L(X_{t^-},\ud y)\,\ud t\biggr],
\end{align}
and $\sfE_\sfR[Z_\infty^v\,|\,\sF_\tau]=Z_\tau^v$ for any $\bF$-stopping time $\tau$.  
\end{proposition}

\subsection{Stochastic optimal control formulations}
\blue
Given a terminal cost $f:\varGamma\to[0,+\infty]$ and initial data $\mu$ satisfying $\Ent(\mu\,|\,\sfR_0)<+\infty$, we consider, as in the introduction, the stochastic optimal control problem
\begin{align}\label{eq:soc}
	\gamma_{\rm sto}(\mu)\coloneqq\inf \left\{\sfE_{\sfQ} \biggl[ f(X_\tau) + \log\frac{\ud \sfQ}{\ud \sfR} \biggr] \; :\; \sfQ\in \sP(\Omega),\;\sfQ_0=\mu \right\}.\tag{{\sf s-OC$_\tau$}}
\end{align}
For a general stopping time $\tau$ and possibly unbounded $f$, defined only on the $\sigma$-algebra $\sigma(X_\tau)$, the following result provides the existence of minimizers to \eqref{eq:soc}.
\black

\begin{theorem}\label{thm:soc}
	Let $(\Omega,\sF,\bF,\sfR)$ be a complete right-continuous filtered space and $\tau$ be an $\bF$-stopping time. Let $f:\varGamma\to[0,+\infty]$ be a proper $\sigma(X_\tau)$-measurable terminal cost satisfying additionally:
	\begin{align}\label{eq:f-property}
		\text{$\exists$ a $\sigma(X_\tau)$-measurable set $D\subset \text{dom}\,f$:}\quad \sfR_\tau(D)>0\quad\text{and}\quad \int_D f\,\ud\sfR_\tau <+\infty.\tag{$\mathsf{A}_f$}
	\end{align}
	Then, for every $x\in\text{dom}\, f$, the stochastic optimal control problem 
\begin{align}\label{eq:soc-x}
	\bar\gamma_{\rm sto}(x)\coloneqq \inf \left\{\sfE_{\sfQ} \biggl[ f(X_\tau) + \log\frac{\ud \sfQ}{\ud \sfR^x} \biggr] \; :\; \sfQ\in \sP(\Omega) \right\}\tag{{\sf s-OC$_\tau^x$}}
\end{align}	
	admits a unique optimal path measure $\sfP^x\in\sP(\Omega)$ with $\sfP_0=\delta_x$.
	
	In particular, the stochastic optimal control problem \eqref{eq:soc} admits a minimizer for any measure $\mu\in\sP(\varGamma)$ with $\Ent(\mu\,|\,\sfR_0)<+\infty$ taking the form $\sfP = \mu\otimes \sfP^x$.
\end{theorem}
\begin{proof}
	  We begin by showing that the variational problem is well-defined, i.e., the set of admissible path measures is nonempty. For this, we consider $D\subset \text{dom}\,f$ as in $\eqref{eq:f-property}$ and set $\sfQ\coloneqq \nu\otimes \sfR^{X_\tau}$, i.e.,
	\[
		\sfQ \coloneqq \int_\varGamma \sfR^{X_\tau=\eta}\,\nu(\ud \eta)\qquad\text{with}\qquad \nu \coloneqq \frac{\mathbf{1}_D}{\sfR_\tau(D)}\sfR_\tau.
	\]
	We then see from \eqref{eq:entropy-disintegrate} with $\phi=X_\tau$ that
	\begin{align*}
		\sfE_{\sfQ} \biggl[ f(X_\tau) + \log\frac{\ud \sfQ}{\ud \sfR} \biggr] &= \sfE_{\sfQ} \biggl[ f(X_\tau) + \log\frac{\ud \sfQ_\tau}{\ud \sfR_{\tau}}(X_\tau)\biggr] = \sfE_\nu\biggl[f + \log \frac{\ud \nu}{\ud \sfR_\tau}\biggr] \\
		&= \int_\varGamma f\,\ud\nu + \Ent(\nu\,|\,\sfR_\tau) = \frac{1}{\sfR_\tau(D)}\int_D f \,\ud\sfR_\tau -\log\sfR_\tau(D)<+\infty,
	\end{align*}
	thus implying that $\sfQ$ is an admissible path measure for the minimization problem \eqref{eq:soc}. 
	
	We now move on to the existence of minimizers for \eqref{eq:soc}. The first part of the statement follows from a straightforward application of the Direct Method of the Calculus of Variations. For completeness, we detail the argument.
	
	Consider a minimizing sequence $(\sfP^k)_{k\ge 1}\subset\sP(\Omega)$ satisfying
	\[
		\sup_{k\ge 1} \sfE_{\sfP^k} \biggl[ f(X_\tau) + \log\frac{\ud \sfP^k}{\ud \sfR^x} \biggr] <+\infty.
	\]
	Since $f\ge 0$, this implies that
	\[
		\sup_{k\ge 1}\Ent(\sfP^k\,|\,\sfR^x) <+\infty.
	\]
	From the de la Vall\'ee-Poussin theorem, we have that the sequence of Radon-Nikodym derivatives
	\[
		\biggl(Z^k\coloneqq\frac{\ud \sfP^k}{\ud \sfR^x}\biggr)_{k\ge 1}\in L^1(\Omega,\sfR^x)\quad\text{is uniformly integrable.}
	\]
	 The Dunford-Pettis theorem then provides the existence of some non-negative  function $Z\in L^1(\Omega,\sfR^x)$ with $\|Z\|_{L^1(\sfR^x)}=1$ and a (not relabelled) subsequence such that
	\[
		Z^k \rightharpoonup Z \quad\text{weakly in $L^1(\Omega,\sfR^x)$}.
	\]
	Setting $\sfP^x\coloneqq Z\sfR^x$, we deduce the setwise convergence $\sfP^k\rightharpoonup \sfP^x$. 
	
	Owing to the lower semicontinuity of the relative entropy w.r.t.\ setwise convergence and the measurability of $\Omega\ni \omega\mapsto f(X_{\tau(\omega)}(\omega))$, we obtain for every $n\ge 1$,
	\begin{align*}
		\sfE_{\sfP^x} \biggl[ (f\wedge n)(X_\tau) + \log\frac{\ud \sfP^x}{\ud \sfR^x} \biggr] &\le \lim_{k\to\infty} \sfE_{\sfP^k} \bigl[ (f\wedge n)(X_\tau)\bigr] + \liminf_{k\to\infty} \Ent(\sfP^k\,|\,\sfR^x) \\
		& \le \liminf_{k\to\infty} \sfE_{\sfP^k} \biggl[ (f\wedge n)(X_\tau) + \log\frac{\ud \sfP^k}{\ud \sfR^x} \biggr] \\
		&\le \liminf_{k\to\infty} \sfE_{\sfP^k} \biggl[ f(X_\tau) + \log\frac{\ud \sfP^k}{\ud \sfR^x} \biggr].
	\end{align*}
	Passing to the limit $n\to\infty$ using the monotone convergence theorem then yields
	\[
		\sfE_{\sfP^x} \biggl[ f(X_\tau) + \log\frac{\ud \sfP^x}{\ud \sfR^x} \biggr] \le \liminf_{k\to\infty} \sfE_{\sfP^k} \biggl[ f(X_\tau) + \log\frac{\ud \sfP^k}{\ud \sfR^x} \biggr],
	\]
	thus implying the minimality of the path measure $\sfP^x\ll \sfR^x$ with $\sfP^x_0=\delta_x$. The strict convexity of the cost functional easily implies the uniqueness of the minimizer $\sfP^x$.

\smallskip

As for the second part, we simply notice that for any $\sfQ\in\sP(\Omega)$ with $\sfQ_0=\mu$ and
\[
	0\le \Ent(\sfQ\,|\,\sfR)=\sfE_{\sfQ} \biggl[ \log\frac{\ud \sfQ}{\ud \sfR} \biggr]\le \sfE_{\sfQ} \biggl[ f(X_\tau) + \log\frac{\ud \sfQ}{\ud \sfR} \biggr].
\]
We have that $\sfQ_0\ll \sfR_0$ and the disintegration theorem gives
\[
	\frac{\ud \sfQ}{\ud \sfR}(\omega) = \frac{\ud \sfQ_0}{\ud \sfR_0}(X_0(\omega))\, \frac{\ud \sfQ^{X_0}}{\ud \sfR^{X_0}}(\omega)\qquad\text{for $\sfQ$-almost every $\omega\in\Omega$.}
\]
Moreover, using formula \eqref{eq:entropy-disintegrate} with $\phi=X_0:\Omega\to \varGamma$, we arrive at the expression
\begin{align*}
	\Ent(\sfQ\,|\,\sfR) =\Ent(\sfQ_0\,|\,\sfR_0) + \int_\varGamma \Ent(\sfQ^x\,|\, \sfR^x)\,\sfQ_0(\ud x).
\end{align*}
Consequently, \eqref{eq:soc} is reduced to
\begin{align*}
  \Ent(\mu\,|\,\sfR_0) + \inf \left\{ \int_\Gamma    \sfE_{\sfQ}^x\biggl[f(X_\tau) + \log \frac{\ud\sfQ}{\ud\sfR^x}  \biggl]\,\mu(\ud x) \right\} . 
\end{align*}
Therefore, by setting $\sfP^x$, $x\in\varGamma$, to be the minimizer of \eqref{eq:soc-x}, we obtain a minimizer of \eqref{eq:soc} by selecting the path measure $\sfP=\mu\otimes \sfP^x$.
\end{proof}

\begin{remark}
	The previous theorem indicates that considering the minimization problem for atomic initial measures $\delta_x$, $x\in\varGamma$, suffices since the optimal path measure for more general admissible initial measures $\mu$ can be obtained from the minimizers of $\eqref{eq:soc-x}$ via superposition.
\end{remark}

\blue
According to Proposition~\ref{prop:girsanov}, \eqref{eq:soc} can be reformulated as a stochastic optimal control problem for controlled path measures $\sfP^v$ with bounded velocity fields $v$, namely
\begin{align}\label{eq:partial-oc-inf-pre}
	\inf \left\{\sfE_{\sfQ} \biggl[ f(X_\tau)+ \iint_{(0,\tau]\times \varGamma}\ent\biggl(\frac{\ud \overline{L}_\sfQ}{\ud \overline{L}_\sfR}\biggr)\, \ud \overline{L}_\sfR\biggr]\; :\; \sfQ\in \sP(\Omega),\; \sfQ_0=\mu \right\}.\tag{{\sf ps-OC$_\infty$}}
\end{align}
Generally, \eqref{eq:partial-oc-inf-pre} has a lower cost than \eqref{eq:soc}. To see this, we use the law of total expectation and Jensen's inequality to obtain
    \begin{align*}
        \sfE_{\sfP^v}\biggl[f(X_\tau) + \log\frac{\ud \sfP^v}{\ud \sfR}\biggr] 
        &= \sfE_{\sfP^v}\bigl[f(X_\tau)\bigr] + \sfE_{\sfR}\biggl[\frac{\ud \sfP^v}{\ud \sfR}\log\frac{\ud \sfP^v}{\ud \sfR}  \biggr]\\
        &= \sfE_{\sfP^v}\bigl[f(X_\tau)\bigr] + \sfE_{\sfR} \Bigl[\sfE_{\sfR}\bigl[\phi(Z_\infty^v)\big| \sF_\tau \bigr] \Bigl],\qquad \phi(s) = s\log s\\
%        &= \sfE_{\sfP^v} \big[f(X_\tau) \big]+ \sfE_{\sfR} \biggl[\sfE_{\sfR}\big[\phi(Z^v) \big| \sF_\tau \big] \biggl]
        &\geq \sfE_{\sfP^v}\bigl[f(X_\tau)\bigr] + \sfE_{\sfR} \Bigl[\phi\big(\sfE_{\sfR}\bigl[ Z_\infty^v  \big| \sF_\tau \bigr]\big) \Bigl] \\
        &= \sfE_{\sfP^v}\bigl[f(X_\tau)\bigr] + \sfE_{\sfR} \bigl[\phi\big(  Z^v_\tau \big) \bigr]
        = \sfE_{\sfP^v}\Bigl[  f(X_\tau)+ \log(Z^v_\tau)\Bigr].
    \end{align*}
From the first conclusion of Proposition \ref{prop:girsanov}, we have that
\[
	\log Z_\tau^v - \iint_{(0,\tau]\times\varGamma} \ent\biggl(\frac{\ud \overline{L}^v}{\ud \overline{L}_\sfR}\biggr)\, \ud \overline{L}_\sfR\quad
%		\underbrace{\sum_{s\in \sT_\tau} \log v_s(X_{s^-},X_s) - \iint_{(0,t]\times\varGamma} \log v_s(X_{s^-},y)\,\overline{L}^v(\ud s\,\ud y)}_
		\text{is a $\sfP^v$-martingale on $\{Z_\tau^v>0\}$} ,
\]
and thus we then obtain the commonly used expression \eqref{eq:partial-oc-inf-pre} for the optimal control formulation with a stopping time.
\black

\section{Application to transition paths: Closed formula solution}\label{sec:sto-tpt}
%\subsection{Application to transition path problems}\label{sec:sto-tpt}
%\begin{remark}\label{rem:ex-transition}

\blue
In this section, we narrow down to the transition path problem which aims to construct an optimal control such that the transition from set $A$ to $B$ can be realized efficiently without altering the bridges. We will follow the running cost derived in Section \ref{sec:sto-opt} and introduce a specific terminal cost. By the additive property of the relative entropy under disintegration, we will verify the bridges for the optimally controlled process does not change and provide a closed formula solution for the optimal path measure $\sfP\!_{AB}$ and optimally controlled L\'evy kernel $L_{AB}$ via committor function $h_{AB}$; see Theorem \ref{thm:sto-rep}.

Precisely, we consider  two nonempty, disjoint and measurable subsets $A$, $B\subset\varGamma$.
	In the context of transition path applications, a penalty is assigned for states in the local attractor $A\subset \varGamma$, while a reward is given in $B\subset\varGamma$, i.e., the desired terminal cost reads
	\begin{equation}\label{tp-f}
    	f_{AB}(x) \coloneqq \begin{cases}
         +\infty & \text{for $x\in A$},  \\
         0 & \text{for $x\in B$}.
    \end{cases}
\end{equation}

Moreover, the values of $f_{AB}$ in $(A\cup B)^c$ should {\em not} influence the cost in any way. This restriction can be implemented as follows: Define the stopping time as the first hitting time $\tau=\tau_{AB}\coloneqq\inf\{t\geq 0 : \, X_t \in  A\cup B\}$ of the set $A\cup B\subset\varGamma$, and consider   the terminal cost $f_{AB}(X_{\tau_{AB}})$. Then only the values of $f_{AB}$ at states in $A$ and $B$ are required.
Therefore, we choose the terminal cost $f=f_{AB}\coloneqq -\log h_{AB}$, where $h_{AB}\ge 0$ is the unique solution of the boundary value problem \eqref{qn}, i.e.,
\begin{equation}\label{eq:elliptic}
	\begin{cases}\displaystyle
		\;\; \int_{\varGamma} \dnabla h(x,y)\,L(x,\ud y) = 0 & \quad x\in  (A\cup B)^c, \\
		\hspace*{7.4em} h(x)=\mathbf{1}_B(x) & \quad x\in A\cup B.
	\end{cases}\tag{{\sf BVP}}
\end{equation}
 
If $(\Omega,\sF,\bF,\sfR)$ is a complete right-continuous filtered space, then the first hitting time $\tau_{AB}$ is an $\bF$-stopping time \cite[Theorem~9.7]{Kallenberg2021}. In this case, $f$ is $\sigma(X_{\tau})$-measurable and Theorem~\ref{thm:soc} provides a solution $\sfP\in\sP(\Omega)$ to the transition path problem formulated as the stochastic optimal control problem \eqref{eq:soc}. 

After choosing the specific stopping time $\tau=\tau_{AB}$ and terminal cost $f=f_{AB}$ in the optimal control formulation \eqref{eq:soc-x}, the disintegration formula leads to an explicit form for the optimal path measure $\sfP\!_{AB}$, which turns out to be the law of the Doob-conditioned process selecting all the paths reaching $B$ before $A$ at $\tau_{AB}$; see \eqref{eq:sto-optimal}. From this point, it is natural to see that the committor function $h_{AB}$ provides the optimal control with associated value function. Thus, we arrive at another main result of this paper, i.e., an explicit representation of the value function 
	\begin{align}\label{eq:soc-gamma}
		\gamma_{\rm sto}(\mu) = \Ent(\mu\,|\,\sfR_0) -\int_\varGamma \log h_{AB}\,\ud \mu,
	\end{align}
This representation of the value function is one of the main motivations of this work as it provides a variational approach in determining the committor function $h_{AB}$. Indeed, if $\sfR_0 = \sfR_0^x = \mu$, $x\in (A\cup B)^c$, then \eqref{eq:soc-gamma} reduces to $\bar\gamma_{\rm sto}(x)=\gamma_{\rm sto}(\delta_x) = -\log h_{AB}(x)$, and therefore,
	\[
		h_{AB}(x) = \exp\bigl(-\bar\gamma_{\rm sto}(x)\bigr),\qquad x\in (A\cup B)^c.
	\]
	%For this reason, we focus on this scenario in the rest of the section.
We now present the formal derivation of the closed-form optimal control. Rigorous proofs are given via passage to the limit in an approximated optimal control problem with cut-off terminal cost $f^n$.
\black
\medskip

The formal idea of obtaining the representation \eqref{eq:soc-gamma} is by means of the disintegration theorem. Indeed, disintegrating any measure $\sfP^x\ll \sfR^x$ along $\sfP_\tau^x$, we obtain
\[
	\frac{\ud \sfP^x}{\ud \sfR^x}(\omega) = \frac{\ud \sfP_\tau^x}{\ud \sfR_\tau^x}(X_\tau(\omega))\,\frac{\ud \sfP^{x,X_\tau}}{\ud \sfR^{x,X_\tau}}(\omega)\qquad \text{for $\sfP^x$-almost every $\omega\in\Omega$}.
\]
Therefore, using \eqref{eq:entropy-disintegrate} with $\phi=X_\tau:\Omega\to\varGamma$, the cost function in \eqref{eq:soc} can be further expressed as
\[
	\sfE_{\sfP^x}\biggl[f(X_{\tau}) + \log\frac{\ud \sfP^x}{\ud \sfR^x}\biggr] = \sfE_{\sfP^x}\biggl[f(X_{\tau}) + \log\frac{\ud \sfP_\tau^x}{\ud \sfR_\tau^x}(X_\tau)\biggr] + \int_\varGamma \Ent(\sfP^{x,X_\tau=\eta}\,|\,\sfR^{x,X_\tau=\eta})\,\sfP_\tau^x(\ud \eta),
\]
where the latter term can be made zero by simply choosing 
\[
	\sfP^{x,X_\tau=\eta}=\sfR^{x,X_\tau=\eta}\qquad\text{for $\sfP_\tau^x$-almost every $\eta\in\varGamma$.}
\]  Moreover, letting
\begin{align}\label{eq:P-tau}
	\frac{\ud \sfP_\tau^x}{\ud \sfR_\tau^x}(X_\tau) = \frac{\exp(-f(X_\tau))}{\sfE_{\sfR^x}[\exp(-f(X_\tau))]} = \frac{h_{AB}(X_\tau)}{\sfE_{\sfR^x}[h_{AB}(X_\tau)]}\qquad\text{$\sfP^x$-almost surely}, 
\end{align}
we then obtain
\[
	\sfE_{\sfP^x}\biggl[f(X_{\tau}) + \log\frac{\ud \sfP^x}{\ud \sfR^x}\biggr] = -\log \sfE_{\sfR^x}\bigl[h_{AB}(X_\tau)\bigr] = - \log h_{AB}(x),
\]
where the last equality follows from Dynkin's formula (cf.\ \eqref{eq:dynkin}), i.e., for $x\in (A\cup B)^c$,
\begin{equation}\label{dynkin}
    \sfE_{\sfR^x}\bigl[h_{AB}(X_\tau)\bigr] - h_{AB}(x) = \sfE_{\sfR^x}\biggl[\iint_{(0,\tau]\times\varGamma} \,\dnabla h_{AB}(X_{s^-},y)\,\overline{L}_\sfR(\ud s\,\ud y)\biggr] = 0.
\end{equation}
From this discussion, we see that a candidate minimizer of \eqref{eq:soc} takes the form
\begin{align}\label{eq:sto-minimizer}
	\sfP^x = \int_\varGamma \sfR^{x,X_\tau=\eta}\,\sfP_\tau^x(\ud \eta) = \frac{h_{AB}(X_\tau)}{\sfE_{\sfR^x}[h_{AB}(X_\tau)]}\sfR^x,
\end{align}
with terminal law $\sfP_\tau^x$ given in \eqref{eq:P-tau}.

\medskip

Our first result  is the explicit formula for the optimal path measure  when the terminal cost is bounded.

\begin{proposition}\label{prop:TPT_n}
	Let $h_{AB}^n$, $n\ge 1$, be the solution of the boundary value problem 
	\begin{equation}\label{eq:elliptic-approx}
	\begin{cases}\displaystyle
		\;\; \int_{\varGamma} \dnabla h(x,y)\,L(x,\ud y) = 0 & \quad x\in  (A\cup B)^c, \\
		\hspace*{7.4em} h(x)=\mathbf{1}_B(x) + e^{-n}\mathbf{1}_A(x) & \quad x\in A\cup B,
	\end{cases}
\end{equation}
	and set $f^n\coloneqq -\log h_{AB}^n$. Then for any $x\in\varGamma$,
\begin{enumerate}[(i)]
\item \eqref{eq:soc-x} admits a unique minimizer given by
\begin{align}\label{eq:sto-optimal-n}
		\sfP_{AB}^{n,x}\coloneqq \frac{h_{AB}^n(X_{\tau_{AB}})}{h_{AB}^n(x)}\sfR^x 
		= \int_\varGamma \sfR^{x,X_{\tau_{AB}}=\eta}\,\frac{h_{AB}^n(\eta)}{h_{AB}^n(x)} \,\sfR_{\tau_{AB}}^x(\ud \eta).
	\end{align}
	\item the associated value function is given by
	\begin{align}\label{eq:soc-gamma-n}
		\bar\gamma^n_{\rm sto}(x) = - \log h^n_{AB}(x).
	\end{align}	
	\item $\sfP_{AB}^{n,x}\in \sP(\Omega)$ solves the martingale problem $\MP(\overline{L}_{AB}^n,\delta_x)$ with L\'evy kernel
\begin{equation}\label{vKernel-n}
	\overline{L}_{AB}^n(\omega,\ud t\,\ud y) \coloneqq \mathbf{1}_{[0,\tau_{AB})}(t)\frac{h_{AB}^n(y)}{h_{AB}^n(X_{t^-}(\omega))} L(X_{t^-}(\omega),\ud y)\,\ud t,\qquad \omega\in \Omega.
\end{equation}
\end{enumerate}
\end{proposition}
\begin{proof}
It is easy to check that property \eqref{eq:f-property} of Theorem~\ref{thm:soc} is satisfied, from which we obtain a unique minimizer for \eqref{eq:soc}. We will use verification procedures to obtain a closed formula for the minimizer and then verify its optimality by disintegrating at the stopping time, which concludes {\em (i)} and {\em (ii)}. Then, we apply the Girsanov transform in Proposition \ref{prop:girsanov} for bounded velocity fields to prove conclusion {\em (iii)}.

\paragraph{\em Step 1}
Let $h_{AB}^n$, $n\ge 1$, be the solution of the boundary value problem \eqref{eq:elliptic-approx}, then the maximum principle for the transition rate $L$ implies $e^{-n}\leq h^n_{AB}\leq 1$. Set $\tau\coloneqq \tau_{AB}$.
	
	For each $n\ge 1$, we set the path measure $\sfP_{AB}^{n,x}$ as in \eqref{eq:sto-optimal-n}. As previously shown,
	\begin{align*}
		 \sfE_{\sfP_{AB}^{n,x}}\biggl[f^n(X_\tau) + \log \frac{\ud\sfP_{AB}^{n,x}}{\ud\sfR^x}\biggr]
		= -\log \sfE_{\sfR^x}\big[h_{AB}^n(X_\tau)\bigr] = -\log h_{AB}^n(x).
	\end{align*}

\paragraph{\em Step 2}
	We now claim that $\sfP_{AB}^{n,x}$ is indeed optimal for each $n\ge 1$. Indeed, let $\sfQ\in\sP(\Omega)$ be any admissible path measure for \eqref{eq:soc-x} with $f^n$ instead of $f$. Setting 
	\begin{align*}
		\sfQ^{\dagger,x} \coloneqq \int_\varGamma \sfR^{x,X_\tau=\eta}\,\sfQ_\tau^x(\ud \eta),
	\end{align*}
	we easily find from \eqref{eq:entropy-disintegrate} that
	\begin{align*}
		\sfE_{\sfQ^x} \biggl[ f^n(X_\tau) + \log\frac{\ud \sfQ^x}{\ud \sfR^x} \biggr] 
		&\ge \sfE_{\sfQ^x} \biggl[ f^n(X_\tau) + \log\frac{\ud \sfQ_\tau^x}{\ud \sfR_\tau^x}(X_\tau) \biggr] \\
		&= \sfE_{\sfQ^{\dagger,x}} \biggl[ f^n(X_\tau) + \log\frac{\ud \sfQ_\tau^{\dagger,x}}{\ud \sfR_\tau^x}(X_\tau) \biggr] + \int_\varGamma \Ent(\sfR^{x,X_\tau}\,|\,\sfR^{x,X_\tau})\,\ud\sfQ_\tau^x \\
		&= \sfE_{\sfQ^{\dagger,x}} \biggl[ f^n(X_\tau) + \log\frac{\ud \sfQ^{\dagger,x}}{\ud \sfR^x} \biggr],
	\end{align*}
	where we used the form of $\sfQ^{\dagger,x}$ and the fact that the relative entropy is non-negative. Furthermore, since $\sfP_{AB,\tau}^{n,x}\sim \sfR_\tau^x$ (i.e., they are mutually absolutely continuous), thus also $\sfQ_\tau^{\dagger,x}\ll \sfP_{AB,\tau}^{n,x}$, and we obtain
	\begin{align*}
	\sfE_{\sfQ^{\dagger,x}} \biggl[ f^n(X_\tau) + \log\frac{\ud \sfQ^{\dagger,x}}{\ud \sfR^x} \biggr]  &= \sfE_{\sfQ^{\dagger,x}} \biggl[ -\log h_{AB}^n(X_\tau) + \log\frac{\ud \sfQ_\tau^{\dagger,x}}{\ud \sfR_\tau^x}(X_\tau) \biggr] \\ &=\int_\varGamma \biggl(-\log \frac{\ud \sfP_{AB,\tau}^{n,x}}{\ud \sfR_\tau^x} + \log\frac{\ud \sfQ_\tau^{\dagger,x}}{\ud \sfR_\tau^x}\biggr)\,\ud \sfQ^{\dagger,x}_\tau - \log \sfE_{\sfR^x}\bigl[h_{AB}^n(X_\tau)\bigr] \\
	&=\int_\varGamma \log \frac{\ud \sfQ_\tau^{\dagger,x}}{\ud \sfP_{AB,\tau}^{n,x}} \,\ud\sfQ^{\dagger,x}_\tau - \log  h_{AB}^n(x) \\
	&\ge - \log h_{AB}^n(x) = \sfE_{\sfP_{AB}^{n,x}}\biggl[f^n(X_\tau) + \log \frac{\ud\sfP_{AB}^{n,x}}{\ud\sfR^x}\biggr],
\end{align*}
where we used \eqref{eq:P-tau} in the second equality.
This verifies the optimality of $\sfP^{n,x}$ for $x\in(A\cup B)^c$, and therewith asserting the representation \eqref{eq:soc-gamma-n}, therewith concluding the proof of statements {\em (i)} and {\em (ii)}.

\paragraph{\em Step 3} Finally, we set the velocity field 
\[
	v_{AB}^n(x,y)\coloneqq  \frac{h_{AB}^n(y)}{h_{AB}^n(x)} \le e^n,\qquad (x,y)\in E_\varGamma,
\]
and apply Girsanov's theorem for bounded velocities (cf.\ Proposition~\ref{prop:girsanov}) to find that
\[
	Z_t^n \coloneqq \exp\left(- \iint_{(0,t]\times \varGamma} \bigl( v_{AB}^n(X_{s^-},y)-1\bigr)\,\overline{L}_\sfR(\ud s\,\ud y)\right)\prod_{s\in \sT_{t}} v_{AB}^n (X_{s^-},X_s) ,
\]
is a strictly positive $\sfR$-martingale w.r.t.\ $\bF$. Since $h_{AB}^n$ solves the boundary value problem \eqref{eq:elliptic}, we further have that for $t\le \tau_{AB}$,
\begin{align*}
	&\iint_{(0,t]\times \varGamma} \bigl(v_{AB}^n(X_{s^-},y)-1\bigr)\,\overline{L}_\sfR(\ud s\,\ud y) \\
	&\hspace{6em}= \iint_{(0,t]\times \varGamma} \frac{\dnabla h_{AB}^n(X_{s^-},y)}{h_{AB}^n(X_{s^-})}\,\overline{L}_\sfR(\ud s\,\ud y) = 0\qquad \text{$\sfR$-almost surely.}
\end{align*}
(Note: the \eqref{eq:elliptic} equation holds on $(A\cup B)^c$, so the above identities are valid for $t\le\tau_{AB}$.) Therefore, $Z_t^n$ simplifies to, for $t\le\tau_{AB}$,
\[
	Z_t^n = \prod_{s\in \sT_{t}} v_{AB}^n (X_{s^-},X_s) = \prod_{s\in\sT_t}\frac{h_{AB}^n(X_s)}{h_{AB}^n(X_{s^-})} = \frac{h_{AB}^n(X_t)}{h_{AB}^n(X_0)}\qquad \text{$\sfR$-almost surely.}
\]
In particular, $\sfP_{AB}^{n,x} = Z_\tau^n\, \sfR^x$ and we deduce from Proposition~\ref{prop:girsanov} that $\sfP_{AB}^{n,x}$ solves the martingale problem $\MP(\overline{L}_{AB}^n,\delta_x)$, thereby concluding the proof.
\end{proof}

We now extend the previous proposition to the unbounded terminal cost $f=f_{AB}$ defined in \eqref{tp-f} via Gamma-convergence.

\begin{theorem}\label{thm:sto-rep}
	Let $h_{AB}$ be the solution of the boundary value problem \eqref{eq:elliptic} and $f=f_{AB}$ be as defined in \eqref{tp-f}. 
	Then for any $x\in (A\cup B)^c$ satisfying $\sfE_{\sfR^x}[e^{\tau_{AB}}]<+\infty$,
\begin{enumerate}[(i)]
\item \eqref{eq:soc-x}  admits a unique  minimizer given by
\begin{align}\label{eq:sto-optimal}
		\sfP_{AB}^x\coloneqq \frac{h_{AB}(X_{\tau_{AB}})}{h_{AB}(x)}\sfR^x
		=\int_\varGamma \sfR^{x,X_{\tau_{AB}}=\eta}\frac{h_{AB}(\eta)}{h_{AB}(x)} \,\sfR_{\tau_{AB}}^x(\ud \eta).
	\end{align}	
	\item the associated value function is given by
	\[
		\bar\gamma_{\rm sto}(x) = - \log h_{AB}(x).
	\]
	\item $\sfP_{AB}^{x}\in \sP(\Omega)$ solves the martingale problem $\MP(\overline{L}_{AB},\delta_x)$ with L\'evy kernel
\begin{equation}\label{vKernel}
	\overline{L}_{AB}(\omega,\ud t\,\ud y) \coloneqq \mathbf{1}_{[0,\tau_{AB})}(t)\frac{h_{AB}(y)}{h_{AB}(X_{t^-}(\omega))} L(X_{t^-}(\omega),\ud y)\,\ud t,\qquad \omega\in \Omega.
\end{equation}
\end{enumerate}	
\end{theorem}
\begin{proof}
	Following the proof of Theorem~\ref{thm:soc}, one can deduce from the existence of some path measure $\sfQ^\circ\in\sP(\Omega)$ such that
	\[
		\sE(\sfQ^\circ)\coloneqq \sfE_{\sfQ^\circ}^x \biggl[ f(X_\tau) + \log\frac{\ud \sfQ^\circ}{\ud \sfR^x} \biggr] <+\infty.
	\]
	We will prove {\em (i)} and {\em (ii)} by means of Gamma-convergence and the results in Proposition~\ref{prop:TPT_n} for the cut-off terminal cost $f^n\coloneqq f \wedge n$. One can then pass to the limit $n\to\infty$ to obtain the required representation. Recall $e^{-n}\leq h_{AB}^n\leq 1$ and $h_{AB}\leq h^{n}_{AB}$.

We define the functionals
\begin{align*}
\sE^n(\sfQ)\coloneqq \sfE_{\sfQ}^x \biggl[ f^n(X_\tau) + \log\frac{\ud \sfQ}{\ud \sfR^x} \biggr],\qquad n\ge 1.
\end{align*}
Clearly, $\sE_n\le \sE$ for every $n\ge 1$. By Proposition \ref{prop:TPT_n}, we have
\begin{equation}\label{tmQQ}
\bar\gamma^n_{\rm sto}(x) = \min_{\sfQ\in\sP(\Omega)} \sE^n(\sfQ) = - \log h^n_{AB}(x) \leq  - \log h_{AB}(x) <+\8.
\end{equation}
We now prove the Gamma-convergence of $\sE^n$ to $\sE$ under setwise convergence and consequently deduce the existence of a minimizer for $\sE$.

\smallskip

\noindent\emph{Liminf inequality:} Let $\sfQ^n \rightharpoonup \sfQ$ setwise in $\sP(\Omega)$. The result is trivial if $\liminf_{n\to\infty} \sE^n(\sfQ^n)=+\8$. Now assume that $\liminf_{n\to\infty} \sE^n(\sfQ^n)<+\8$ and set $m\ge 1$. Since $\sfE_{\sfQ^n}^x[f^n(X_\tau)] \geq \sfE_{\sfQ^n}^x [f^m(X_\tau)]$ for every $n\ge m$, we have that
\[
	\liminf_{n\to\infty} \sfE_{\sfQ^n}^x [f^n(X_\tau)] \geq \liminf_{n\to\infty} \sfE_{\sfQ^n}^x [f^m(X_\tau)] = \sfE_\sfQ^x [f^m(X_\tau)],
\]
where the last equality follows from the setwise convergence. By the monotone convergence theorem, we then have that
\begin{equation}\label{finf}
	\liminf_{n\to\infty} \sfE_{\sfQ^n}^x [f^n(X_\tau)] \geq \sfE_\sfQ^x [f(X_\tau)].
\end{equation}
Owing to the lower semicontinuity of relative entropy, we conclude with \eqref{finf} that
\[
	\sE(\sfQ)\le \liminf_{n\to\infty} \sE^n(\sfQ^n). 
\]

\smallskip

\noindent\emph{Limsup inequality:} Let $\sfQ\in \sP(\Omega)$. The limsup inequality is trivial if $\sE(\sfQ)=+\8$. Now suppose $\sE(\sfQ)<+\8$. For the constant sequence $\sfQ^n\coloneqq \sfQ$, we have by construction $\sE^n(\sfQ^n)\le \sE(\sfQ)$. Hence,
\[
	\limsup_{n\to\infty} \sE^n(\sfQ^n) \leq \sE(\sfQ).
\]

\smallskip

\noindent{\em Existence of minimizer:} For each $n\ge 1$, we let $\sfP_{AB}^{n,x}$ be the unique minimizer of $\sE^n$ given in Proposition~\ref{prop:TPT_n}. Since $f^n\ge 0$, and 
\[
	\sup_{n\ge 1}\Ent(\sfP_{AB}^{n,x}\,|\,\sfR^x)\le \sup_{n\ge 1}\sE^n(\sfP_{AB}^{n,x}) \le \sup_{n\ge 1}\sE^n(\sfQ^\circ) \le \sE(\sfQ^\circ) <+\infty,
\]
we have that the sequence $(Z_\tau^n = \ud\sfP_{AB}^{n,x}/\ud \sfR^x)_{n\ge 1}\subset L^1(\Omega,\sfR^x)$ is uniformly integrable. The Dunford-Pettis theorem then provides the existence of some non-negative function $Z_\tau\in L^1(\Omega,\sfR^x)$ with $\|Z_\tau\|_{L^1(\sfR^x)}=1$ and a (not relabelled) subsequence such that
\[
	Z_\tau^n\rightharpoonup Z_\tau\quad\text{weakly in $L^1(\Omega,\sfR^x)$}\qquad \text{with}\qquad \int_\Omega \ent(Z_\tau)\,\ud\sfR^x<+\infty,
\]
where we recall $\ent(s) = s\log s -s +1$.
Setting $\sfP_{AB}^x\coloneqq Z_\tau\,\sfR^x$, we deduce the setwise convergence $\sfP_{AB}^{n,x}\rightharpoonup \sfP_{AB}^x$. The Gamma-convergence result then implies that $\sfP_{AB}^x$ is a minimizer of $\sE$, and since $\sE$ is strictly convex, it is the unique minimizer. Moreover, the 
%Since $f^n\geq 0$, the boundedness of the relative entropy implies there exists a (not relabelled) subsequence $\sfQ^n$ setwise in $\sP(\Omega)$. Thus the sub-level set of $\sE^n$ is pre-compact in the sense of the setwise convergence in $\sP(\Omega)$.
%Therefore, we conclude the minimizer  $P^n_{AB} = \argmin \sE^n(Q)$ converges to the minimizer  $P = \argmin \sE(Q)$. Notice the minimizer of $\sE(Q)$ is also unique due to strictly convexity of the relative entropy. 
associated value function for $\sE^n$ also converges to the minimum value of $\sE$, i.e.,
\[
	\bar\gamma_{\rm sto}(x)=\lim_{n\to\infty} \bar\gamma^n_{\rm sto}(x) = - \lim_{n\to\infty} \log h^n_{AB}(x) = - \log h_{AB}(x),
\]
where we used the $h_{AB}^n\to h_{AB}$ pointwise, thus concluding statement {\em (ii)}. 

From the pointwise convergence $h_{AB}^n\to h_{AB}$, we further obtain 
\[
	Z_\tau^n = \frac{h_{AB}^n(X_\tau)}{h_{AB}^n(x)} \longrightarrow  \frac{h_{AB}(X_\tau)}{h_{AB}(x)}=Z_\tau\qquad\text{$\sfR^x$-almost surely}.
\]
where the right-most equality holds due to the uniqueness of weak limits, yielding {\em (i)}.

\smallskip

Finally, we prove point {\em (iii)}.
Define
 \[
	v_{AB}(x,y)\coloneqq  \frac{h_{AB}(y)}{h_{AB}(x)},\qquad (x,y)\in E_\varGamma.
\] The property $\sfP_{AB}^x(X_0\in\cdot)=\delta_x$ directly follows by construction. Since $\overline L_{AB}$ is possibly unbounded, the proof of the second property will require some preparation. We set
\begin{align*}
	\Delta_k\coloneqq \inf \Bigl\{t>0: \overline L_{AB}((0,t\wedge \tau]\times\varGamma)>2^k\Bigr\},\qquad k\ge 1.
\end{align*}
Since $\overline L_{AB}$ is a predictable Markov L\'evy kernel, $(\Delta_k)_{k\ge 1}$ is a sequence of stopping times, which we claim to converge to infinity $\sfP_{AB}^x$-almost surely. Indeed, from $\sfE_{\sfR^x}[e^\tau]<+\infty$, we obtain
\begin{align*}
	\sfE_{\sfP_{AB}^x}[\tau] = \sfE_{\sfR^x}[\tau Z_\tau] \le \sfE_{\sfR^x}[\ent^*(\tau)] + \Ent(\sfP_{AB}^x\,|\,\sfR^x) <+\infty,
\end{align*}
where we made use of Fenchel-Young inequality for the conjugate pair $(\ent,\ent^*)$, with $\ent^*(s)=e^s-1$. For each $k\ge 1$, we then estimate
\begin{align*}
	\sfP_{AB}^x(\Delta_k<+\infty) &= \sfP_{AB}^x\bigl(\overline L_{AB}((0,\Delta_k\wedge\tau]\times\varGamma)>2^k\bigr) \\
	&\le \frac{1}{2^k}\sfE_{\sfP_{AB}^x}\bigl[\,\overline L_{AB}((0,\Delta_k\wedge\tau]\times\varGamma)\bigr] \\
	&\le \frac{1}{2^k}\sfE_{\sfP_{AB}^x}\biggl[ c_L\ent^*(1)\,\tau + \iint_{(0,\tau]\times\varGamma}\ent(v_{AB}(X_{s^-},y))\, L(X_{s^-}(\omega),\ud y)\,\ud s\biggr] \\
	&= \frac{1}{2^k}\Bigl( c_L\ent^*(1)\,\sfE_{\sfP_{AB}^x}[\tau] + \Ent(\sfP_{AB}^x\,|\,\sfR^x)\Bigr) <+\infty,
\end{align*}
where we again used the Fenchel-Young inequality for $(\ent,\ent^*)$. Hence, we deduce
\[
	\sum_{k=1}^\infty \sfP_{AB}^x(\Delta_k<+\infty) \le c_L\ent^*(1)\,\sfE_{\sfP_{AB}^x}[\tau] + \Ent(\sfP_{AB}^x\,|\,\sfR^x)<+\infty.
\]
By the Borel-Cantelli Lemma, $\lim_{k\to\infty}\Delta_k =+\infty$ $\sfP_{AB}^x$-almost surely as desired.
	
\smallskip
	
Now let $\varphi\in B_b(\bR^+\times E_\varGamma)$ be a non-negative bounded Borel function and set
\begin{align*}
	N_t^{\varphi}\coloneqq \sum_{s\in \sT_{t}}\varphi_s(X_{s^-},X_s) - \iint_{(0,t]\times \varGamma} \varphi_s(X_{s^-},y)\,\overline{L}_{AB}(\ud s\,\ud y),\qquad t\ge 0.
\end{align*}
Owing to the boundedness of the process $(N_{t\wedge\Delta_k}^\varphi)_{t\ge 0}$, we can follow the proof of Proposition~\ref{prop:girsanov} to infer that $(N_{t\wedge \Delta_k}^\varphi Z_{t\wedge\tau})_{t\ge 0}$ is an $\sfR^x$-martingale w.r.t.\ $\bF$, where $Z_{t\wedge\tau}=\sfE_{\sfR^x} [ Z_\tau |  \sF_{t\wedge \tau}]$. Using the fact that
\[
	\iint_{(0,t\wedge\Delta_k]\times \varGamma} \varphi_s(X_{s^-},y)\,\overline{L}_{AB}(\ud s\,\ud y)\, \longrightarrow \iint_{(0,t]\times \varGamma} \varphi_s(X_{s^-},y)\,\overline{L}_{AB}(\ud s\,\ud y)
\]
$\sfP_{AB}^x$-almost surely, and the estimate
\begin{align*}
	&\iint_{(0,t\wedge\Delta_k]\times \varGamma} \varphi_r(X_{r^-},y)\,\overline{L}_{AB}(\ud r\,\ud y) \le c_L\ent^*\bigl(\|\varphi\|_\infty\bigr)
	\,\tau  \\
	&\hspace{8em}+ \iint_{(0,\tau]\times\varGamma}\ent(v_{AB}(X_{s^-},y))\, L(X_{s^-}(\omega),\ud y)\,\ud s\qquad \text{$\sfP_{AB}^x$-almost surely,}
\end{align*}
where the right-hand side is known to be $\sfP_{AB}^x$-integrable, we then apply the dominated convergence theorem to conclude that $(N_t^\varphi)_{t\ge 0}$ is a $\sfP_{AB}^x$-martingale w.r.t.\ $\bF$, therewith concluding the proof.
\end{proof}

\blue
\begin{proof}[Proof of Theorem~\ref{mainthm:tpt}] For notational simplicity, we denote $\tau=\tau_{AB}$. Let $\sfQ$ be an admissible path measure for \eqref{eq:soc-intro}, then \eqref{eq:entropy-disintegrate} gives
\[
    \Ent(\sfQ\,|\,\sfR) = \Ent(\mu\,|\,\sfR_0) + \int_\varGamma \Ent(\sfQ^x\,|\,\sfR^x)\,\mu(\ud x).
\]
In particular, the cost for the variational problem \eqref{eq:soc-intro} may be expressed as
\begin{align}\label{mainthm:tpt:1}
    \sfE_{\sfQ} \biggl[ f(X_\tau) + \log\frac{\ud \sfQ}{\ud \sfR} \biggr] = \Ent(\mu\,|\,\sfR_0) + \int_\varGamma \sfE_{\sfQ}^x \biggl[ f(X_\tau) + \log\frac{\ud \sfQ^x}{\ud \sfR^x} \biggr]\,\mu(\ud x).
\end{align}
Since the integrand in the second term is non-negative, a minimizer of \eqref{eq:soc-intro} is attained by minimizing \eqref{eq:soc-x} for $\mu$-almost every $x\in\varGamma$. According to Theorem~\ref{thm:sto-rep}(i), the unique minimizer takes the form
\[
    \sfP_{AB}^x = \frac{h_{AB}(X_{\tau_{AB}})}{h_{AB}(x)}\sfR^x\qquad\text{for $\mu$-almost every $x\in\varGamma$}.
\]
Setting
\[
    \frac{\ud \sfP\!_{AB}}{\ud \sfR} = \frac{\ud \mu}{\ud \sfR_0}(X_0)\frac{h_{AB}(X_{\tau})}{h_{AB}(X_0)},
\]
then yields a posteriori the unique minimizer of \eqref{eq:soc-intro}. Together with Theorem~\ref{thm:sto-rep}(ii), \eqref{mainthm:tpt:1} yields the explicit form of the value function
\[
    \gamma_{\rm sto}(\mu) = \Ent(\mu\,|\,\sfR_0) -\int_\varGamma \log h_{AB}\,\ud \mu.
\]
Finally, the form of the L\'evy kernel is precisely the one given in Theorem~\ref{thm:sto-rep}(iii).
\end{proof}
\black

\begin{remark}\label{rem_sig}
Notice that the basic Girsanov transform in Proposition~\ref{prop:girsanov} requires the boundedness of the velocity field $v$. However, the transition path problem produces an optimal control velocity field $v_{AB}$ that is unbounded. Hence, Proposition~\ref{prop:girsanov} cannot be directly applied. Yet, in statement (iii) of Theorem \ref{thm:sto-rep}, we proved that the optimal change of measure formula is given by
\begin{equation*}
\sfP_{AB}^x\coloneqq Z_\tau\,\sfR^x, \qquad Z_\tau=\frac{h_{AB}(X_\tau)}{h_{AB}(x)}.
\end{equation*}
Due to the uniform integrability of $Z$, the process defined by
\[
	Z_t \coloneqq \sfE_{\sfR^x} [ Z_\tau |  \sF_{t\wedge \tau}],\qquad t\ge 0,
\]
then gives a non-negative $\sfR$-martingale, which can be used to show that $\sfP_{AB}^x$ solves the martingale problem $\MP(\overline{L}_{AB},\delta_x)$. In a similar spirit, one could eventually extend Proposition~\ref{prop:girsanov} to include a class of singular velocity fields (see \cite{leonard2012girsanov} for a general representation formula for path measures having finite relative entropy).
\end{remark}

\begin{remark}
	We emphasize that the proof of Theorem~\ref{thm:sto-rep} illustrates that the minimization problem \eqref{eq:soc-x} on path measures can be {\em greatly} reduced to solving the minimization problem
	\begin{align}\label{eq:soc-reduced}
		-\log h_{AB}(x)=\inf \left\{\sfE_{\nu} \biggl[ f + \log\frac{\ud \nu}{\ud \sfR_{\tau_{AB}}^x} \biggr] \; :\; \nu\in \sP(\varGamma)\right\}.
	\end{align}
	Indeed, a minimizer $\nu_{AB}$ of \eqref{eq:soc-reduced} defines to a path measure via
	\[
		\sfQ = \int_\varGamma \sfR^{x,X_\tau=\eta}\,\nu_{AB}(\ud \eta),
	\]
	which provides a minimizer of \eqref{eq:soc-x}. The right-hand side of \eqref{eq:soc-reduced} might also provide an alternative approach to determining the value of the committor function in point $x$ without the need to solve the boundary value problem \eqref{eq:elliptic}. 
\end{remark}

\blue
\begin{remark}
In the finite time horizon case, where $\tau=T\in (0,+\infty)$, same as the infinite time horizon optimal control formulation \eqref{eq:partial-oc-inf-pre}, we have the finite time horizon optimal control problem
\begin{align}\label{eq:partial-oc}
	\inf \left\{\sfE_{\sfQ} \biggl[ f(X_T)+ \iint_{(0,T]\times \varGamma}\ent\biggl(\frac{\ud \overline{L}_\sfQ}{\ud \overline{L}_\sfR}\biggr)\, \ud \overline{L}_\sfR\biggr] \; :\; \sfQ\in \sP(\Omega),\;\sfQ_0=\mu \right\},\tag{{\sf ps-OC$_T$}}
\end{align}
where $\overline{L}_\sfQ$ is the L\'evy kernel associated to $\sfQ$.
Furthermore, suppose the {\em velocity field} $v:(0,T)\times E_\varGamma \to [0,\infty)$   defines a L\'evy kernel 
\[
	\overline{L}^v(\omega,\ud t\,\ud y) \coloneqq \,v_t(X_{t^-}(\omega),y) L(X_{t^-}(\omega),\ud y)\,\ud t,\qquad \omega\in \Omega.
\]
One can also use the density-flux pair $(p,q)$ associated with the forward Kolmogorov equation to convert 
\eqref{eq:partial-oc} to a deterministic optimal control problem. Suppose $v_t \coloneqq \ud q_t/\ud p_t\otimes L$ for almost every $t\in(0,T)$, then we have the {\em deterministic} optimization problem with an ODE constraint
\begin{equation}\label{oc}
   \gamma_{\rm det}(\mu) \coloneqq \inf  \left\{\,\int_{\varGamma} f\, \ud p_T + \int_0^T \Ent(q_t\,|\,p_t\otimes L)\,\ud t\; :\; (p,q)\in \sA_{[0,T]}(\mu)\right\}, \tag{{\sf OC$_{T}$}}
\end{equation}
where $f:\varGamma\to[0,+\infty]$ is a lower semicontinuous {\em terminal cost}. Here, $\sA_{[0,T]}(\mu)$ are admissible density-flux pairs defined in \eqref{aa}. In Appendix \ref{sec:deterministic}, we present existence, uniqueness results, and a closed-form solution via a Hamilton-Jacobi equation \eqref{tHJE} with the terminal condition $\psi_T = - f$. Proofs will also be provided for the reader's convenience.
\end{remark}
\black

\section{Conclusion and Further Discussions}\label{sec5}

We formulate the transition path problem for general jump processes as a stochastic optimal control problem by choosing the relative entropy of path measures as running cost and choosing an infinite penalty as a terminal cost at a stopping time. 
Theorem \ref{thm:sto-rep} provides an explicit solution for this optimal control problem in the infinite time horizon scenario and offers a closed-form expression for the optimal control via the discrete committor function $h_{AB}$. This result permits direct usage of the optimal path measure $\sfP\!_{AB}$ and the corresponding optimally controlled L\'evy kernel $\overline L_{AB}$ for Monte Carlo simulations of transition paths.

Our optimal control formulation, supported by Theorem \ref{thm:sto-rep}, establishes a theoretical foundation for simulating transition paths in various scientific problems modeled by jump processes. This approach is similar to the optimal control interpretation of the committor function for continuous transition paths modeled by diffusion processes. The optimal control formulation for drift-diffusion processes at infinite time horizons has been explored in \cite{GLLL23}. Using the optimally controlled process, the authors proposed a finite temperature mean transition path algorithm, applying it to the folding computation of the Alanine dipeptide. In the zero-diffusion limit, the transition path simulated from the optimally controlled drift-diffusion process recovers the most probable path (minimum action path) computed from the Freidlin-Wentzell large deviation rate function for the exit problem.

For both continuous processes and jump processes, one can observe that the optimally controlled generator is given by the Doob $h$-transformation of the original generator. This is essentially because the equivalence between the change of measures on path spaces and the change in the L\'evy kernel for the controlled process is established via the Girsanov transform for both continuous and pure jump processes. 

If the L\'evy kernel is induced from a spatial discretization of the Fokker-Planck equation for drift-diffusion processes, see for instance \cite{GLLL23}, which uses Voronoi tessellation and corresponding graph derivatives to construct a continuous time Markov chain on Voronoi cell centers, then the optimally controlled processes can be constructed via the discrete committor function $h_{AB}$ via Theorem \ref{thm:sto-rep}. As the Voronoi tessellation approaches a refined limit, the discrete committor function will converge to the continuous one. An intriguing question is whether (and in which sense) the value function and associated transition path will converge in the diffusion limit; see \cite{hraivoronska2022diffusive} for similar results in the context of generalized gradient flows.

We also comment on the relation between the finite time horizon deterministic optimal control problem and the stochastic optimal control problem with a stopping time. Based on our choice of running cost, both of these optimal control problems fit into the same framework via a convex optimization problem for path measures. This relies on the disintegration formula and its additive property for the relative entropy. The problem on the convergence from the finite time horizon deterministic optimal control problem to the infinite time horizon one involves the long-time behavior of the HJE with a terminal value, which we leave for future study. We refer to \cite{namah1999remarks, ichihara2009long, cgmt15} for long time behavior of HJEs and refer to \cite{GL23w} for the selection principle of the stationary solution to HJEs resulting from the exponential change of variable in the linear backward equation for chemical reactions.

Various other methodologies exist for transition path computations, such as the minimum action method \cite{MAM}, the (finite temperature) string method \cite{string, FTSM}, and the transition path theory grounded in the committor function \cite{weinan2006, MMS2009, Lu_Nolen_2015}. In most of those methods, including our optimally controlled Monte Carlo simulations, the committor function computation serves as a vital initial step. Recent computational methods for high dimensional committor functions are explored in \cite{lai2018point, khoo2019solving, li2019computing, chen2023committor}. In Monte Carlo simulations, the optimal change of measures is also known as importance sampling, cf. \cite{BD_book}. Given an inaccurate committor function, how to design an important sampling that stimulates the Monte Carlo simulations for transition paths and updates simultaneously the committor function could be an interesting future study.

\appendix

\section{Girsanov's Theorem: Proof of Proposition~\ref{prop:girsanov}}\label{app:girsanov}

Part of the proof follows the strategy outlined in \cite{leonard2012girsanov} (see also \cite{Jacod1975}). As in Section~\ref{sec:jump-process}, $\sT\coloneqq\{t\ge 0: X_t\ne X_{t^-}\}$ denotes the set of canonical jump times and $\sT_t\coloneqq \sT\cap [0,t]$ the set of jump times restricted to the interval $[0,t]$ for each $t\ge 0$, with $\sT_0=\emptyset$. For any path measure $\sfP\in\sP(\Omega)$ and its corresponding L\'evy kernel $\overline{L}_\sfP$, the random measure
\[
	\overline{\Lambda}_\sfP(\omega,\ud t\,\ud y) \coloneqq \delta_{X_t(\omega)}(\ud y)\,\delta_{\sT(\omega)}(\ud t) - \overline{L}_\sfP(\omega,\ud t\,\ud y),\quad \omega\in \Omega.
\]
denotes the compensated sum of jumps, where the first term on the right-hand side is the canonical random measure \cite{Jacod1975,Kallenberg2021}.

To show that $Z^v$ is a non-negative $\sfR$-martingale, we apply It\^o's formula on $\exp(Y_t^v)$, with
\[
	Y_t^v \coloneqq \sum_{s\in \sT_{t}} \log v_s(X_{s^-},X_s) - \iint_{(0,t]\times \varGamma} \bigl(v_s(X_{s^-},y)-1\bigr)\,\overline{L}_\sfR(\ud s\,\ud y)\qquad\text{ on $\{Z_t^v>0\}$}.
\]
Notice the first term is the pure jump part and the second term is the continuous part with bounded variation, so we deduce the stochastic differential equation
\begin{align*}
	\ud\exp(Y_t^v) &= \Bigl(\exp(Y_{t}^v)-\exp(Y_{t^-}^v)\Bigr)\delta_\sT(\ud t) - \exp(Y_{t^-}^v)\int_\varGamma\bigl(v_t(X_{t^-},y)-1\bigr)\,L(X_{t^-},\ud y)\,\ud t \\
	&= \exp(Y_{t^-}^v) \int_{\varGamma} \bigl( v_t(X_{t^-},y)-1\bigr)\, \overline{\Lambda}_\sfR(\ud t\,\ud y).
\end{align*}
Consequently, we find that $Z^v=\exp(Y^v)$ solves the stochastic equation
\begin{align*}\label{eq:exponential-sde}
	Z_t^v = 1 + \iint_{(0,t]\times \varGamma} Z_{s^-}^v\bigl(v_s(X_{s^-},y)-1\bigr)\, \overline{\Lambda}_\sfR(\ud s\,\ud y).
\end{align*}
Since the integrand of the $\sfR$-stochastic integral above is a bounded predictable process, the second term on the right is an $\sfR$-martingale, thus concluding that $Z^v$ is an $\sfR$-martingale w.r.t.\ $\bF$. Finally, the non-negativity follows easily from construction.

\medskip

Next, we show that $\sfP^v$ defined as in Proposition~\ref{prop:girsanov} solves the martingale problem $\MP(\overline{L}^v,\mu)$. By construction, the first property of Definition~\ref{def:martingale} is satisfied. As for the second property, we set for any bounded measurable function $\varphi\in B_b(\bR^+\times E_\varGamma)$,
\begin{align*}
	N_t^\varphi \coloneqq
	\sum_{s\in \sT_{t}}\varphi_s(X_{s^-},X_s) - \iint_{(0,t]\times \varGamma} \varphi_s(X_{s^-},y)\,\overline{L}^v(\ud s\,\ud y),
\end{align*}
and show that the product $N^\varphi Z^v$ is an $\sfR$-martingale w.r.t.\ $\bF$. To do so, we make use of the It\^o's formula again for the product to deduce
\begin{align*}
	\ud \bigl(N_t^\varphi Z_t^v\bigr) 
	&= Z_{t^-}^v \int_\varGamma \varphi_t(X_{t^-},y)\,v_t(X_{t^-},y)\,\overline{
	\Lambda}_\sfR(\ud t\,\ud y) + N_{t^-}^\varphi Z_{t^-}^v \int_\varGamma \bigl(v_t(X_{t^-},y)-1\bigr)\,\overline{\Lambda}_\sfR(\ud t\,\ud y),
\end{align*}
which then implies
\[
	N_t^\varphi Z_t^v = \iint_{(0,t]\times\varGamma} Z_{s^-}^v\Bigl[\varphi_s(X_{s^-},y)\,v_s(X_{s^-},y) + N_{s^-}^\varphi\bigl(v_s(X_{s^-},y)-1\bigr)\Bigr]\,\overline{
	\Lambda}_\sfR(\ud s\,\ud y).
\]
Since the integrand on the right-hand side is a bounded predictable process, the right-hand side is an $\sfR$-martingale w.r.t.\ $\bF$, and thus also the left-hand side. We then conclude that the process $N^\varphi$ is a $\sfP^v$-martingale w.r.t.\ $\bF$ as desired.

\medskip

Finally, suppose $\Ent(\sfP^v\,|\,\sfR)<+\infty$. To show the formula for the relative entropy, we use \eqref{eq:entropy-disintegrate} with $\phi=X_0$ and \eqref{eq:P-girsanov-formula} to obtain
\[
	\Ent(\sfP^v\,|\,\sfR) = \Ent(\mu\,|\,\sfR_0) + \sfE_{\sfR} \bigl[Z^v\log Z^v\bigr].
\]
Since both terms on the right-hand side are non-negative, they are both finite due to the finiteness of the left-hand side. In particular, the $\sfR$-martingale $Z^v$ is uniformly integrable in $L^1(\Omega,\sfR)$, and Doob's martingale convergence theorem \cite[Theorem~9.22]{Kallenberg2021} provides a limit $Z_\infty^v\in L^1(\Omega,\sfR)$ such that $Z_t^v\to Z_\infty^v$ $\sfR$-almost surely and in $L^1(\Omega,\sfR) $ as $t\to +\infty$. The representation $\sfE_\sfR[Z_\infty^v\,|\,\sF_\tau]=Z_\tau^v$ follows from the optional sampling theorem \cite[Theorem~7.18]{Fima_book}.

Since for every $t\ge 0$,
\[
	\log Z_t^v 
		= \underbrace{\sum_{s\in \sT_t} \log v_s(X_{s^-},X_s) - \iint_{(0,t]\times\varGamma} \log v_s(X_{s^-},y)\,\overline{L}^v(\ud s\,\ud y)}_{\sfP^v\text{-martingale on $\{Z_t^v>0\}$}} + \iint_{(0,t]\times\varGamma} \ent\biggl(\frac{\ud \overline{L}^v}{\ud \overline{L}_\sfR}\biggr)\, \ud \overline{L}_\sfR,
\]
we then obtain the desired expression for the relative entropy.  \qed

\section{Deterministic Optimal Control Formulation}\label{sec:deterministic}

\blue In this section, we give results and proofs for the  existence and closed formula solution of the  deterministic optimal control problem \eqref{oc}.  We will obtain the explicit formula for the optimal control and the value function can be obtained by directly solving a linear backward equation, which is obtained via the Cole-Hopf transformation for the nonlinear  Hamilton-Jacobi equation.
Particularly, in Proposition \ref{prop:finite-unbounded}, we provide a rigorous justification for the existence of minimizers and the formula for the optimal value function via Gamma-convergence when the terminal cost takes infinite values.
\black

%\subsection{Deterministic formulation and Hamilton--Jacobi equation}\label{sec_HJE}
We first present an optimal control problem defined by an ODE-constrained optimization problem in a fixed time horizon, where the ODE is motivated by the evolution of the time-marginal law $\sfP_t\coloneqq(X_t)_\#\sfP$ of a path measure $\sfP\in \sP(\Omega)$ given by the {\em forward Kolmogorov equation}
\begin{align}\label{eq:FKE}
	\int_\varGamma \varphi(x)\,\sfP_t(\ud x) -\int_\varGamma \varphi(x)\,\sfP_s(\ud x) = \iint_{(s,t)\times \varGamma\times\varGamma} \dnabla\varphi(x,y)\, L(x,\ud y)\sfP_r(\ud x)\ud r,\tag{{\sf FKE}}
\end{align}
for every $\varphi\in B_b(\varGamma)$ and interval $[s,t]\subset[0,T]$.

\smallskip

In view of \eqref{eq:FKE}, we make the following definition.

\begin{definition}\label{def:CE}
	A density-flux pair $(p,q)$, where
	\begin{enumerate}
		\item $p\in C([0,T];\sP(\varGamma))$ is continuous curve w.r.t.\ setwise convergence, and
		\item $q\in L^1((0,T);\sM_+(E_\varGamma))$, $E_\varGamma\coloneqq\varGamma\times\varGamma$, is an integrable family of fluxes,
	\end{enumerate}
	is said to satisfy the discrete {\em continuity equation}
	\begin{align}\label{eq:CE}
		\partial_t p_t + \ddiv q_t = 0,\qquad \ddiv q(\ud x)\coloneqq \int_{y\in \varGamma} \bigl[q(\ud x\,\ud y) - q(\ud y\,\ud x)\bigr],\tag{{\sf CE}}
	\end{align}
	if for every $[s,t]\subset[0,T]$ and $\varphi\in B_b(\varGamma),$
	\[
		\int_\varGamma \varphi\,\ud p_t - \int_\varGamma \varphi\,\ud p_s = \int_s^t\!\! \int_{E_\varGamma} \dnabla \varphi(x,y)\,q_r(\ud x\,\ud y)\,\ud r.
	\]
	
	A density-flux pair $(p,q)$ is said to have {\em finite action} if, additionally,
\begin{align}\label{eq:action}
	\mathscr{L}_{[0,T]}(p,q) \coloneqq \int_0^T \Ent(q_t\,|\,p_t\otimes L)\,\ud t <+\infty,\qquad p\otimes L(\ud x\ud y) \coloneqq L(x,\ud y)\,p(\ud x).
\end{align}
We call a density-flux pair $(p,q)$ {\em admissible} if it satisfies the discrete continuity equation \eqref{eq:CE} and has finite action \eqref{eq:action}; we denote by $\sA_{[0,T]}(\mu)$ the family of all admissible density-flux pairs with initial density $p_0=\mu\in\sP(\varGamma)$, i.e.,
\begin{align}\label{aa}
\begin{aligned}
\sA_{[0,T]}(\mu) &\coloneqq\Bigl\{ (p,q) \;\text{satisfying \eqref{eq:CE}, \eqref{eq:action}, and $p_0=\mu$}
	 \Bigr\}\\
	 &\hspace{10em} \subset C([0,T];\sP(\varGamma))\times \sM_+((0,T)\times E_\varGamma).
\end{aligned}
\end{align}
Furthermore, we say that a sequence of admissible pairs $\{(p^n,q^n)\}_{n\ge 1}\subset \sA_{[0,T]}(\mu)$ converges to  an admissible pair $(p,q)\in \sA_{[0,T]}(\mu)$ if 
\begin{gather}\label{eq:limits} \left\{\quad
			\begin{gathered}
			p_t^n \rightharpoonup p_t\quad\text{setwise in $\sP(\varGamma)$ for every $t\in[0,T]$},\\
			q^n\rightharpoonup q\quad \text{setwise in $\sM((0,T)\times E_\varGamma)$}.
			\end{gathered}\right.
		\end{gather}
		We denote this convergence simply with $(p^n,q^n)\rightharpoonup (p,q)$ in $\sA_{[0,T]}(\mu)$.
\end{definition}

\begin{remark}
	If $\Upsilon$ is countable, then setwise convergence in $\sM(\Upsilon)$ coincides with the strong convergence on $\ell^1(\Upsilon)$. In particular, functions $f:[0,T] \to \ell^1(\Upsilon)$ are strongly continuous if and only if they are setwise continuous.
\end{remark}

% \begin{remark}\label{rem:admissible-levy}
% We will see later that any admissible density-flux pair $(p,q)$ defines a stochastic process. Indeed, since $\mathscr{L}_{[0,T]}(p,q)<+\infty$, the {\em velocity field} $v:(0,T)\times E_\varGamma \to [0,\infty)$, with $v_t \coloneqq \ud q_t/\ud p_t\otimes L$ for almost every $t\in(0,T)$, defines a L\'evy kernel 
% \[
% 	\overline{L}^v(\omega,\ud t\,\ud y) \coloneqq \,v_t(X_{t^-}(\omega),y) L(X_{t^-}(\omega),\ud y)\,\ud t,\qquad \omega\in \Omega.
% \]
% A solution $\sfP$ of the martingale problem associated with $\overline{L}^v$ then gives rise to a jump process.  
% \end{remark}

An important property of the class $\sA_{[0,T]}$ of admissible density-flux pairs is the following compactness and lower semicontinuity result, whose proof is an adaptation of \cite[Proposition~4.21]{PRST2022}.

\begin{lemma}\label{lem:compact-action}
	Let $\{(p^n,q^n)\}_{n\ge 1}\subset \sA_{[0,T]}(\mu)$ with $\mu\ll \pi$ be a family of admissible pairs:
	\[
		\sup_{n\ge 1}\mathscr{L}_{[0,T]}(p^n,q^n) <+\infty.
	\]
	Then there exists a (not relabelled) subsequence and an  admissible pair $(p,q)\in \sA_{[0,T]}(\mu)$ such that
\[
	(p^n,q^n)\rightharpoonup (p,q)\quad\text{in $\sA_{[0,T]}(\mu)$ in the sense of \eqref{eq:limits}.}
\]

		Additionally, for any sequence of admissible pairs with $(p^n,q^n)\rightharpoonup(p,q)$ in $\sA_{[0,T]}(\mu)$, we have
		\[
			\mathscr{L}_{[0,T]}(p,q)\le \liminf_{n\to\infty}\mathscr{L}_{[0,T]}(p^n,q^n),
		\]
		i.e., the action functional $\mathscr{L}_{[0,T]}$ is lower semicontinuous under the convergence \eqref{eq:limits}.
\end{lemma}

\medskip

The following proposition states that the optimal control for the deterministic problem \eqref{oc} attains a minimizer, which follows from the Direct Method of the Calculus of Variations.

\begin{proposition}\label{prop:existence-oc}
	Let $f:\varGamma\to[0,+\infty]$ be a measurable terminal cost \blue with $\text{dom}\,f\ne \emptyset$. \black Then, for each $T<+\infty$, the finite-time-horizon problem \eqref{oc} admits an optimal admissible pair $(p,q)\in \sA_{[0,T]}(\mu)$. 
\end{proposition}

In Proposition \ref{thm_finite_pde} below, we show that the optimal control flux field $q$ can be described by the solutions of a Hamilton-Jacobi equation and the feedback control takes a time-dependent  Doob-transformation form \eqref{Q_hn}. Before that, let us first introduce essential equations.

Consider the {\em backward Kolmogorov equation}
\begin{equation}\label{eq:bke}
	\partial_t h_t(x)+\int_{\varGamma}\dnabla h_t(x,y)\, L(x,\ud y)=0,\qquad x\in\varGamma,\tag{{\sf BKE}}
\end{equation}
with terminal condition $h_t|_{t=T} = h_T$. Recall that $\dnabla\varphi(x,y) = \varphi(y)-\varphi(x)$. As in the case of the boundary value problem \eqref{qn}, we will assume that a unique classical solution to \eqref{eq:bke} exists and that the maximum principle and instantaneous positivity hold for \eqref{eq:bke}. In particular, we assume that if $0\le h_T\le M_h$ for some constant $M_h> 0$ and $h_T\not\equiv 0$, then $h_t\in (0,M_h]$ for all $t\in[0,T)$.

We introduce an exponential change of variable $\psi_t=\log h_t$, known also as the Cole-Hopf transform or Varadhan's nonlinear semigroup,
which solves the dynamical Hamilton-Jacobi equation for the potential field $\psi$ in the classical sense:
\begin{equation}\label{tHJE}
 \partial_t\psi_t(x) +  H(x, \dnabla\psi_t)=0,\qquad (t,x)\in(0,T)\times\varGamma, \qquad \lim_{t\to T_-} \psi_t = \log h_T.\tag{{\sf HJE}}
\end{equation} 
Here $H:\varGamma\times B_b(E_\varGamma)\to [0,+\infty)$ is the {\em local Hamiltonian} defined by
\[
	H(x, \xi) \coloneqq \int_\varGamma \bigl(\exp(\xi(x,y))-1\bigr)L(x,\ud y),\qquad x\in\varGamma,
\]
which is convex in the $\xi$-argument for every $x\in\varGamma$. 

We further define the {\em global Hamiltonian} $\sH:\sP(\varGamma)\times B_b(E_\varGamma)\to [0,+\infty)$ as
\[
	\sH(p,\xi) \coloneqq \int_{\varGamma} H(x,\xi)\,p(\ud x) = \int_{E_\varGamma} \bigl(e^\xi-1\bigr) \,\ud p\otimes L.
\]
The global Hamiltonian gives rise to a {\em Lagrangian} $\sL:\sP(\varGamma)\times \sM_+(E_\varGamma)\to[0,+\infty]$ via Legendre-Fenchel duality,
\begin{align}\label{tmL}
	\sL(p,q) &\coloneqq \sup_{\xi\in B_b(E_\varGamma)} \left\{ \int_{E_\varGamma} \xi\,\ud q -  \sH(p,\xi) \right\} 
	= \begin{cases}\displaystyle
		\int_{E_\varGamma}\ent\left(\frac{\ud q}{\ud p\otimes L}\right)\ud p\otimes L &\text{if $q\ll p\otimes L$}, \\
		+\infty &\text{otherwise,}
	\end{cases}
\end{align}
where the supremum is achieved at $\xi = \log (\ud q/\ud p\otimes L)$ for  $q\ll p\otimes L$.
Notice that the rightmost expression in \eqref{tmL} is precisely the integrand of the action in \eqref{eq:action}, i.e.,
\[
	\sL(p,q) = \Ent(q\,|\,p\otimes L),
\]
therewith motivating the running cost $\mathscr{L}_{[0,T]}$ in \eqref{oc}.

Utilizing the Hopf-Cole transform and reducing the solution of \eqref{tHJE} to the solution of \eqref{eq:bke} allows us to obtain an explicit formula for the optimal control whenever the terminal cost $f$ is a bounded measurable function.
%{\red Below mimic Proposition \ref{prop:TPT_n}} 
\begin{proposition}\label{thm_finite_pde} Let $f:\varGamma\to [0,\infty)$ be a bounded measurable terminal cost. Then,
\begin{enumerate}[(i)]
	\item there exists a unique classical solution $\psi:[0,T]\times\varGamma\to \bR$ to the Hamilton-Jacobi equation \eqref{tHJE} with terminal condition $\psi_T = - f$.
	\item there is a unique curve $p\in C^1([0,T];\sP(\varGamma))$ satisfying the Cauchy problem
	\begin{align}\label{eq:HJE-flux}
		\partial_t p_t + \ddiv  q_t = 0,\qquad  q_t = \exp(\dnabla\psi_t)\, p_t\otimes L,
	\end{align}
	with initial condition $ p_{t}|_{t=0} = \mu\in\sP(\varGamma)$.
	\item the density-flux pair $( p, q)\in \sA_{[0,T]}(\mu)$ is the unique minimizer  for \eqref{oc}.
	\item the associated value function can be expressed as
	\begin{equation}\label{gamma}
		\gamma_{\rm det}(\mu) = - \int_\varGamma \psi_0\,\ud \mu.
	\end{equation}
\end{enumerate}
\end{proposition}
\begin{proof}
	Let $h:[0,T]\times\varGamma\to \bR$ be the unique solution of the backward Kolmogorov equation \eqref{eq:bke} with terminal condition $h_T=e^{-f}$. Since $0\le f\le M_f\coloneqq\sup_{x\in \varGamma} f(x)$, we have that $e^{-M_f}\le h_T\le 1$. The maximum principle for \eqref{eq:bke} then yields $e^{-M_f}\le h_t\le 1$ for every $t\in[0,T]$. Hence, the function $-M_f\le \psi_t \coloneqq\log h_t \le 0$ satisfies \eqref{tHJE} in the pointwise sense, as derived previously. Using this $\psi_t$, the existence and the uniqueness of the curve $p_t$ solving the Cauchy problem \eqref{eq:HJE-flux} is obvious. Thus, we obtain (i)-(ii).

Next, by verification, we will prove the optimality of the density-flux pair $( p, q)\in \sA_{[0,T]}(\mu)$ in (iii) and the associated value function in (iv).

	Let $(\mu,\eta)\in \sA_{[0,T]}(\mu)$ be an arbitrary admissible density-flux pair.
	Multiplying the Hamilton-Jacobi equation \eqref{tHJE} with $\mu_t$ and integrating over $t$ yields
	\begin{align*}
		\int_0^T\!\!\!\int_\varGamma\partial_t \psi_t\,\ud \mu_t\,\ud t + \int_0^T \sH(\mu_t,\dnabla\psi_t)\,\ud t = 0.
	\end{align*}
	Integrating by parts and using the continuity equation \eqref{eq:CE} then gives
	\begin{align}\label{eq:CE-HJE}
		 - \int_\varGamma \psi_0\,\ud \mu = \int_\varGamma f\,\ud \mu_T +\int_0^T \left\{\int_{E_\varGamma} \dnabla\psi_t\,\ud \eta_t -\sH(\mu_t,\dnabla\psi_t)\right\}\ud t.
	\end{align}
	Using the dual characterization of $\sL$, we then obtain
	\[
		- \int_\varGamma \psi_0\,\ud \mu \le \int_\varGamma f\,\ud \mu_T +\int_0^T \sL(\mu_t,\eta_t)\,\ud t.
	\]
	Since the left-hand side is independent of $(\mu,\eta)$, taking the infimum over density-flux pairs yields
	\[
		- \int_\varGamma \psi_0\,\ud \mu \le \inf_{(\mu,\eta)} \left\{\int_\varGamma f\,\ud \mu_T +\int_0^T \sL(\mu_t,\eta_t)\,\ud t\right\} =\gamma_{\rm det}(\mu).
	\]
	
	Now, let $p\in C([0,T];\sP(\varGamma))$ be a solution of \eqref{eq:HJE-flux}, which exists since the time-dependent jump kernel is uniformly bounded, i.e.
	\[
		\sup_{t\in[0,T]}\sup_{x\in\varGamma}\int_\varGamma \exp(\dnabla\psi_t(x,y))\, L(x,\ud y) < +\infty.
	\]
	Moreover, it is clear that $(p,q)$ with $q_t = \exp(\dnabla\psi_t)\, p_t\otimes L$ has finite action, with $\xi = \dnabla\psi $ achieving the supremum in \eqref{tmL}. Inserting this pair $(p, q)$ in \eqref{eq:CE-HJE}, we obtain
	\[
		- \int_\varGamma \psi_0\,\ud \mu = \int_\varGamma f\,\ud p_T +\int_0^T \left\{\int_{E_\varGamma} \dnabla\psi_t\,\ud q_t -\sH( p_t,\dnabla\psi_t)\right\}\ud t 
		= \int_\varGamma f\,\ud p_T + \mathscr{L}_{[0,T]}( p, q),
	\]
	i.e., the density-flux pair $( p, q)$ is optimal and \eqref{gamma} holds.
\end{proof}

\begin{remark}  The \blue optimal control $v$ used to regulate the transition rate $L$ is defined as $v_t(x,y) \coloneqq \frac{\ud q_t}{\ud p_t\otimes L}(x,y)$ for almost every $t\in(0,T).$ Then, as an essential consequence of Proposition~\ref{thm_finite_pde},
\begin{align}\label{eq:velocity-representation}
	v_t(x,y) \coloneqq \frac{\ud q_t}{\ud p_t\otimes L}(x,y) = e^{\psi_t(y) - \psi_t(x)} = \frac{h_t(y)}{h_t(x)},\quad (x,y)\in E_\varGamma,
\end{align}\black
where $\psi$ solves \eqref{tHJE} and $h$ solves \eqref{eq:bke} with terminal cost $\psi_T=-f$ and $h_T=e^{-f}$, respectively. This enables us to restrict the control variable to a smaller class in the minimization problem \eqref{oc} and the controlled transition rate follows the Doob $h$-transformation formula. 
More specifically, we can express the control variable in terms of the variable $h$ (or equivalently the potential $\psi$). With this control, \blue the controlled transition rate becomes
\begin{equation}\label{Q_hn}
     L^h_t(x,\ud y) \coloneqq \frac{h_t(y)}{h_t(x)} L(x,\ud y),\qquad x\in\varGamma.
\end{equation}
This is the counterpart of the last statement in Theorem \ref{thm:sto-rep}. \black
The corresponding running cost for the density-flux pair $(p,q)$ with $q_t=p_t\otimes L_t^h$, in terms of the $h$-variable, takes the form
\begin{equation}\label{L_hn}
	\mathscr{L}_{[0,T]}(p,q) = \int_0^T\!\!\!\int_{E_\varGamma} \ent\left(\frac{h_t(y)}{h_t(x)}\right) L(x,\ud y)\,p_t(\ud x)\,\ud t. 
\end{equation}
\end{remark}

Next, we consider the case when $f$ is unbounded but proper, i.e., $f\not\equiv +\8$. From the instantaneous positivity property, the Cauchy problem \eqref{eq:HJE-flux} is still well-posed when the terminal cost $f=-\log h_T\not\equiv +\8$. However, it is not clear whether the control velocity taking the form \eqref{eq:velocity-representation} remains optimal. Below, we use the Gamma-convergence technique to give a positive answer. The representation \eqref{gamma} of the value function $\gamma_{\rm det}(\mu)$ still holds, as shown in the following statement.

\begin{proposition}\label{prop:finite-unbounded}
	Let $f:\varGamma\to[0,+\infty]$ be a proper, measurable terminal cost. Then,
\begin{enumerate}[(i)]
	\item there exists a unique classical solution $\psi:[0,T)\times\varGamma\to \bR$ to the Hamilton-Jacobi equation \eqref{tHJE} with terminal condition $\lim_{t\to T^-}\psi_t = - f$;
	\item there exists a unique optimal density-flux pair $( p, q)\in \sA_{[0,T]}(\mu)$ which solves \eqref{oc}.
	\item this curve $p\in C^1([0,T);\sP(\varGamma))$ satisfies the Cauchy problem
	\begin{align}\label{eq:HJE-flux-n}
		\partial_t p_t + \ddiv  q_t = 0,\quad \text{ with }\,  q_t = \exp(\dnabla\psi_t)\, p_t\otimes L,
	\end{align}
	and initial condition $ p_{t}|_{t=0} = \mu\in\sP(\varGamma)$;
	\item the associated value function can be expressed as
	\begin{equation}\label{gamma-n}
		\gamma_{\rm det}(\mu) = - \int_\varGamma \psi_0\,\ud \mu.
	\end{equation}
\end{enumerate}
\end{proposition}
\begin{proof}
We prove the statement via a cut-off approximation and Gamma-convergence argument of the corresponding `cut-off' functionals.
	The existence of the minimizer is a consequence of the Gamma-convergence, while its uniqueness follows from the strict convexity of the cost functional. The explicit formula for the optimal density-flux pair \eqref{eq:HJE-flux-n} will be derived via the limit of \eqref{eq:HJE-flux}.
	
	\paragraph{\em Step 1} Let $f^n\coloneqq f\wedge n$.  Let $h$, $h^n$, $n\ge 1$ be the solution to the backward Kolmogorov equation \eqref{eq:bke} with terminal condition $h_T=e^{-f}$ and $h_T^n=e^{-f^n} = h_T \vee e^{-n}$, respectively. Due to the linearity of \eqref{eq:bke} and the maximum principle, we easily establish that $h_t^n(x)\le h_t^{n+1}(x) \le h_t(x)$, $n\ge 1$, for all $x\in\varGamma$ and $t\in[0,T]$, i.e., for each $t\in[0,T]$, the sequence $(h_t^n)_{n\ge 1}$ is a pointwise non-decreasing sequence.  Hence, the function $\psi_t^n\coloneqq\log h_t^n$ converges pointwise and monotonically to $\psi_t\coloneqq\log h_t$, in which $\psi_t^n$ and $\psi_t$ satisfies \eqref{tHJE} with terminal conditions $\psi_T^n=-f^n$ and $\lim_{t\to T^-}\psi_t=-f$, respectively. This completes the proof of statement {\em (i)}.
	
	In addition, since $L$ is bounded operator in $B_b(\varGamma)$ with bound $c_L$, we have that
	\begin{equation}\label{strong_h}
	\|h_t - h_t^n\|_{\sup} = \|e^{(T-t)L}(h_T - h_T^n)\|_{\sup} \leq e^{(T-t)c_L}\|h_T - h_T^n\|_{\sup} \leq e^{-n} e^{(T-t)c_L},
\end{equation}	
i.e., $h^n$ converges to $h$ uniformly in $[0,T]\times \varGamma $.
	
	\paragraph{\em Step 2}
	
	For each $n\ge 1$, we then obtain from Proposition~\ref{thm_finite_pde} the expression
		\begin{equation}\label{vg1}
		\gamma_{\rm det}^n(\mu) = -\int_\varGamma\psi_0^n\,\ud \mu,
		\end{equation}		 
		where $\gamma_{\rm det}^n(\mu)$ is the value function for the optimal control problem with terminal cost $f^n$. The monotone convergence theorem then implies
		\begin{equation}\label{vg2}
		\lim_{n\to\infty} \int_\varGamma (-\psi_0^n)\,\ud \mu = \int_\varGamma(-\psi_0)\,\ud \mu <+\infty,
		\end{equation}
		where we used the instantaneous positivity property for \eqref{eq:bke}, giving $\psi_0 =\log h_0>-\infty$.
	
	\paragraph{\em Step 3} In the following, we define the functionals 
	\[
		\sA_{[0,T]}(\mu)\ni (p,q)\mapsto \left\{\quad
		\begin{aligned}
			\mathcal{E}^n(p,q) &\coloneqq \int_\varGamma f^n\,\ud p_T + \mathscr{L}_{[0,T]}(p,q),\quad n\ge 1, \\
			\mathcal{E}(p,q) &\coloneqq \int_\varGamma f\,\ud p_T + \mathscr{L}_{[0,T]}(p,q),
		\end{aligned}\right.
	\]
	and show that the sequence of functionals $(\mathcal{E}^n)_{n\ge 1}$ Gamma-converges to $\mathcal{E}$ and that $(\mathcal{E}^n)_{n\ge 1}$ is equi-coercive in the class of admissible pairs $\sA_{[0,T]}(\mu)$. A consequence of this is the convergence of the minimizers and of minima. In particular, $\gamma_{\rm det}^n(\mu)\to \gamma_{\rm det}(\mu)$, which then concludes the proof of the statement {\em (iv)}.

\smallskip

\noindent\emph{Liminf inequality:} Let $(p^n,q^n)\rightharpoonup (p,q)$ in $\sA_{[0,T]}(\mu)$. The lower semicontinuity of $\mathscr{L}_{[0,T]}$ is directly provided by Lemma~\ref{lem:compact-action}. As for the other term, we fix some $m\ge 1$. For every $n\ge m$, we have
		\[
			\int_\varGamma f^n\,\ud p_T^n = \int_\varGamma f^m\,\ud p_T^n + \int_\varGamma (f^n-f^m)\,\ud p_T^n \ge \int_\varGamma f^m\,\ud p_T^n.
		\]
		Since $f$ is measurable, $f^m$ is bounded and measurable, $m\ge 1$, and hence,
		\[
			\liminf_{n\to\infty}\int_\varGamma f^n\,\ud p_T^n \ge \int_\varGamma f^m\, \ud p_T.
		\]
		An application of the monotone convergence theorem for $m\to 0$, then yields
		\[
			\liminf_{n\to\infty}\int_\varGamma f^n\,\ud p_T^n \ge \int_\varGamma f\, \ud p_T.
		\]
		Finally, we obtain
		\[
			\liminf_{n\to\infty} \mathcal{E}^n(p^n,q^n) \ge \liminf_{n\to\infty}\int_\varGamma f^n\,\ud p_T^n + \liminf_{n\to\infty} \mathscr{L}_{[0,T]}(p^n,q^n) \ge \mathcal{E}(p,q).
		\]
		Notice that for the special sequence $(p^n,q^n)$ achieving the value function $\gamma_{\rm det}^n(\mu)$, we have from \eqref{vg2} that  $\lim_{n\to \8} \gamma_{\rm det}^n(\mu) <+\8.$ Therefore, we know $\liminf_{n\to\infty} \mathcal{E}^n(p^n,q^n)<+\8$ and thus
		$$\mathcal{E}(p,q)\leq \liminf_{n\to\infty} \mathcal{E}^n(p^n,q^n)<+\8.$$
		This implies that the set $\{\mathcal{E}(p,q)<+\8\}$ is non-empty and $\gamma_{\rm det}(\mu)<+\8.$
		
\smallskip
		
\noindent\emph{Limsup inequality:} Let $(p,q)$ be an admissible pair. If $\mathcal{E}(p,q)=+\infty$, then there is nothing to prove. If $\mathcal{E}(p,q)<+\infty$, then clearly,
		\[
			\mathcal{E}^n(p,q) \le \mathcal{E}(p,q)\quad\text{for every $n\ge 1$},
		\]
		and hence the limsup inequality is trivially satisfied. Together, the liminf and limsup inequalities provide the Gamma-convergence of $(\mathcal{E}^n)_{n\ge 1}$ to $\mathcal{E}$.% w.r.t.\ the convergence in $\sA_{[0,T]}(\mu)$.
		
\smallskip
		
\noindent\emph{Equi-coercivity:} Since $f, f^n\ge 0$, $n\ge 1$, sublevel sets of $\mathcal{E}^n$ are sublevel sets of $\mathscr{L}_{[0,T]}$ in the class of admissible density-flux pairs, which is non-empty pre-compact due to Lemma~\ref{lem:compact-action}. As a consequence of the Gamma-convergence and equi-coercivity, we conclude that the minimizer $(p^n,q^n)$ of $\mathcal{E}^n$ converges to $(p,q)$ in $\sA_{[0,T]}(\mu)$ in the sense of \eqref{eq:limits} and $(p,q)\in \sA_{[0,T]}(\mu)$ is the minimizer of $\mathcal{E}(p,q)$, therewith concluding the proof of statement {\em (ii)}.

\paragraph{\em Step 4} Recall the explicit formula for the minimizer $(p^n,q^n)$ of $\mathcal{E}^n$ given in \eqref{eq:HJE-flux}, 
\begin{equation}\label{tm_hq}
h^n_t(x) q^n_t(\ud x \ud y) = h_t^n(y) L(x,\ud y) p_t^n(\ud x).
\end{equation}  
From the uniform convergence \eqref{strong_h} and the setwise convergence of $(p^n,q^n)$ to $(p,q)$ in $\sA_{[0,T]}(\mu)$, we have the setwise convergence of both sides in \eqref{tm_hq} so that
\[
h_t(x) q_t(\ud x \ud y) = h_t(y) L(x,\ud y) p_t(\ud x).
\]
Due to the instantaneous positivity of $h_t,\, t<T$, we conclude the explicit formula
\[
	q_t(\ud x \ud y) = \frac{h_t(y)}{h_t(x)} L(x,\ud y) p_t(\ud x)\qquad\text{for $t<T$}.
\]
Recall also that $(p,q)\in \sA_{[0,T]}(\mu)$. Hence, assertion {\em (iii)} is proved.
\end{proof}

\begin{remark}\label{rem1}
Notice that the value function \eqref{gamma} exhibits a linear dependence on the initial distribution $\mu\in\sP(\varGamma)$. If one sets $\mu = \delta_x$, $x\in\varGamma$, then the associated value function becomes
\[
	\bar\gamma_{\rm det}(x)\coloneqq \gamma_{\rm det}(\delta_x) = -\phi_0(x),\qquad z\in\varGamma.
\]
Hence, without loss of generality, one can consider the optimization problem with the jump process conditioned to start from a point $x\in \varGamma$ and compute $\bar\gamma_{\rm det}(x)$. The value function for a general initial distribution $p_0$ can then be readily obtained from $\gamma_{\rm det}(\mu) = \int_\varGamma\bar\gamma_{\rm det}(x)\,\mu(\ud x)$.
\end{remark}

\medskip

In the finite time horizon setting, the controlled generator is expressed in the Doob-transformed form as shown in equation \eqref{Q_hn}, along with the associated running cost presented in entropy form in \eqref{L_hn}. These derivations can be found in \cite[Sec VI.9]{fleming06}; see also \cite{Sheu85}. However, it’s important to note that the terminal cost associated with transition path problems is generally only defined on a specific subset of $\varGamma$ via a stopping time and not everywhere. This limitation restricts the applicability of the optimal control formulation outlined in equation \eqref{oc}.

\section*{Acknowledgments}
Yuan Gao was partially supported by NSF under award DMS-2204288 and NSF CAREER award DMS-2440651. Jian-Guo Liu was partially supported by NSF under award DMS-2106988. Oliver Tse received support from NWO Vidi grant 016.Vidi.189.102 on {\em Dynamical-Variational Transport Costs and Application to Variational Evolutions}, and NWO grants OCENW.M.21.012 and NGF.1582.22.009. Oliver Tse is thankful to Christian L\'eonard for many stimulating discussions.

%%%%%%%%%%%%%%%%%%%%%%%%%%%%%%%%%%%%%%%%%%%%%%%%%%%%%%%%%%%%%%%%%%%
%%                                                               %%
%% You may add acknowledgments (optional).                       %%
%%                                                               %%
%%%%%%%%%%%%%%%%%%%%%%%%%%%%%%%%%%%%%%%%%%%%%%%%%%%%%%%%%%%%%%%%%%%

%%%%%%%%%%%%%%%%%%%%%%%%%%%%%%%%%%%%%%%%%%%%%%%%%%%%%%%%%%%%%%%%%%%
%%                                                               %%
%% You have reached the end of your document.                    %%
%%                                                               %%
%%%%%%%%%%%%%%%%%%%%%%%%%%%%%%%%%%%%%%%%%%%%%%%%%%%%%%%%%%%%%%%%%%%

\end{document}